\documentclass[11 pt]{article}
\usepackage{amsfonts,amssymb}
\usepackage{lipsum,amsmath,amsthm}
\usepackage{ dsfont }
\usepackage{appendix}
\usepackage{color}
\usepackage{esint}
\usepackage{mathtools}
\usepackage{scalerel}
\usepackage{comment}
\usepackage[right = 2.5cm, left=2.5cm, top = 2.5cm, bottom =2.5cm]{geometry}
\usepackage[utf8]{inputenc}
\usepackage[english]{babel}

\usepackage{graphicx} % Required for inserting images

\usepackage[utf8]{inputenc}
\usepackage[english]{babel}

\usepackage{color}
\usepackage{enumerate}
\usepackage{enumitem}
\usepackage{fancyhdr}
\usepackage{upgreek}
\usepackage[T1]{fontenc}

\usepackage{accents}
\usepackage{centernot}
\usepackage{stmaryrd}
\usepackage[super]{nth}

\usepackage{amsmath}
\usepackage{amscd}
\usepackage{amsfonts}
\usepackage{amsthm}
\usepackage{amstext}
\usepackage{amssymb}
\usepackage{amsbsy}
\usepackage{mathrsfs}
\usepackage{mathtools}
\usepackage{array}
\usepackage{tikz-cd}
\usepackage{bm}
 \usepackage{url}

\def\XXint#1#2#3{{\setbox0=\hbox{$#1{#2#3}{\int}$ }
\vcenter{\hbox{$#2#3$ }}\kern-.6\wd0}}

\usepackage{float}
\usepackage{caption}
\usepackage{subcaption}
\usepackage{graphicx}

\usepackage{tikz}
\usetikzlibrary{decorations.markings}
\usetikzlibrary{decorations.pathreplacing}
\usetikzlibrary{calligraphy}
\usetikzlibrary{positioning}
\usetikzlibrary{arrows}
\usetikzlibrary{fadings}
\usetikzlibrary{math}

\usepackage[nobysame]{amsrefs}

\usepackage[linktocpage]{hyperref}
\hypersetup{
	colorlinks=false,
	linkcolor=blue,
	citecolor=red,
	filecolor=cyan,    
	urlcolor=magenta,
}

\makeatletter
\renewcommand*{\@fnsymbol}[1]{\ensuremath{\ifcase#1\or \or \or \or 
   *\or \dagger\or \ddagger\or \mathsection \else\@ctrerr\fi}}
\makeatother

\newtheorem{theorem}[subsection]{Theorem}

\newtheorem{proposition}[subsection]{Proposition}
\newtheorem{corollary}[subsection]{Corollary}

\theoremstyle{definition}

\theoremstyle{remark}
\newtheorem{remark}[subsection]{Remark}

\begin{document}
\title{Mixing Estimates for Passive Scalar Transport by $BV$ Vector Fields\footnote{\textit{Date}: \today.}} % Quantitative well-posedness for passive scalar transport by $BV$ vector fields
\author{Lucas Huysmans\thanks{Max Planck Institute for Mathematics in the Sciences, Inselstra{\ss}e 22, 04103 Leipzig, Germany. \hfill \break \textit{Email address}: \texttt{lucas.huysmans@mis.mpg.de}.} \and Ayman Rimah Said\thanks{Chargé de recherche au  CNRS at Laboratoire de Mathématiques de Reims (LMR - UMR 9008). \hfill \break \textit{Email address}: \texttt{ayman.said@cnrs.fr}.}}
\date{ \ }
\maketitle

\begin{abstract}
    We prove a quantitative mixing estimate for the Cauchy problem for transport along divergence-free vector fields with bounded variation. By developing a framework that quantifies Ambrosio's regularisation scheme, we derive the first explicit bounds on the mixing rate for general $BV$ vector fields. Our analysis reveals that tetration (repeated exponentiation) emerges in the mixing rate from the local nature of Ambrosio's regularisation.
    
    % We give a quantitative mixing/weak stability estimate for passive scalar transport by a time-independent divergence-free vector field with bounded variation. The corresponding bound on the mixing rate is tetration in the time variable and the $BV$-norm of the vector field.
\end{abstract}

%\setcounter{tocdepth}{1}{\small\tableofcontents}

%\newpage
\section{Introduction}
    In this paper, we revisit the Cauchy problem for transport along a divergence-free vector field with bounded variation. It is a remarkable result, due to Ambrosio \cite{ambrosio2004transport}, that these assumptions on the vector field are sufficient for well-posedness. A natural next step is to control the permitted growth of solutions. Since transport conserves the mass of solutions, growth should be measured in the weak topology, e.g., in negative Sobolev norms. Such estimates are more commonly known as mixing estimates since they bound the transfer of mass from high to low frequencies, i.e. how quickly the passive scalar (un)mixes.

    Such mixing estimates are well-studied for Sobolev vector fields in $W^{1,p}$ with $p>1$. However, when $p=1$, or more generally for $BV$ vector fields, there is a significant gap in applying the same harmonic arguments, and no general mixing bounds are known. This problem is further related to an unsolved conjecture by Bressan \cite{bressan2003lemma} on the maximal mixing rate for $W^{1,1}$ vector fields.
    
    In this paper, we develop a novel framework to bridge this gap. By quantifying Ambrosio's well-posedness proof with new harmonic estimates in $BV$ and employing a combinatorial argument with the pigeonhole principle to make these estimates small, we give the first quantitative mixing estimates for passive scalar transport by $BV$ vector fields.

    More precisely, consider a weak solution $\rho(x,t):\mathbb{T}^d\times[0,T]\to\mathbb{R}$ on the $d$-dimension torus $\mathbb{T}^d$, with $\rho(x,t)\in L_t^\infty L_x^\infty$, to the transport equation
    \begin{equation}\label{eq:TE}
        \frac{\partial \rho}{\partial t}(x,t)+\nabla \cdot(u(x)\rho(x,t))=0,
    \end{equation}
    along a divergence-free, time-independent vector field $u(x):\mathbb{T}^d\to\mathbb{R}^d$, with $u(x)\in BV$. Denote by $\rho_0(x)\in L^\infty$ the initial datum, and also by $\rho_T(x)\in L^\infty$ the weak trace of $\rho(x,t)$ at time $t=T$ (see for instance \cite{ambrosio2014continuity}).

    Our primary result is a well-posedness/weak stability estimate which quantifies the continuous dependence of the solution on the initial datum in the negative Sobolev space $W^{-1,1}$:
    {\renewcommand{\thesubsection}{A}
    \begin{theorem}\label{thm:quantitativeBV}
        Let $\kappa>0$. For constants $A,B>0$ depending only on the dimension $d\ge2$, if
        \begin{equation*}
            \left\|\rho_0\right\|_{W^{-1,1}} \le \kappa \left(\exp^{\left\lceil A\exp\left(\frac{B}{\kappa}\left\|\rho_0\right\|_{L^\infty}T\left\|\nabla u\right\|_{L^1}\right)\right\rceil}\left(\frac{1}{\kappa}\left\|\rho_0\right\|_{ L^\infty}\right)\right)^{-1},
        \end{equation*}
        where $\exp^n(x)$ refers to repeated exponentiation $\underbrace{\exp(\exp(\dots(\exp(x))))}_{n}$ (tetration), then
        \begin{equation*}
            \left\|\rho_T\right\|_{W^{-1,1}} \le \kappa.
        \end{equation*}
    \end{theorem}
    }

    By reversing time, one may rephrase Theorem \ref{thm:quantitativeBV} as a more classical bound on mixing, that is, as a \textit{lower} bound on the weak norm $\left\|\rho_T\right\|_{W^{-1,1}}$:
    {\renewcommand{\thesubsection}{B}
    \begin{corollary}[Mixing bound]\label{thm:mixingBV}
        For constants $A,B>0$ depending only on the dimension $d\ge2$, then
        \begin{equation*}
            \left\|\rho_T\right\|_{W^{-1,1}} \ge \left\|\rho_0\right\|_{W^{-1,1}} \left(\exp^{\left\lceil A\exp\left(B\frac{\left\|\rho_0\right\|_{L^\infty}}{\left\|\rho_0\right\|_{W^{-1,1}}}T\left\|\nabla u\right\|_{L^1}\right)\right\rceil}\left(\frac{\left\|\rho_0\right\|_{ L^\infty}}{\left\|\rho_0\right\|_{W^{-1,1}}}\right)\right)^{-1},
        \end{equation*}
        where $\exp^n(x)$ refers to repeated exponentiation $\underbrace{\exp(\exp(\dots(\exp(x))))}_{n}$ (tetration).
    \end{corollary}
    }
    
    The tetration in the bound comes from a pigeonhole principle argument and is related to the \textit{local} in space nature of Ambrosio's regularisation scheme, see Remark \ref{remark:doubleexponential}.

    In contrast, when the vector field lies in the Sobolev space $W^{1,p}$ with $p>1$, one can give exponential mixing bounds by more classical harmonic techniques, see for instance \cite{crippa2008estimates,seis2013maximal,iyer2014lower,leger2018new,hadvzic2018singular,huysmans2024quantitative}. We also mention that in the $BV$ case, one has exponential bounds when restricting the vector field to alternating shears on $\mathbb{T}^2$ \cites{cooperman2023exponential}, or to cellular $90$ degree rotations \cite{hadvzic2018singular}.

    Other work on transport by $BV$ vector fields includes the extension of the well-posedness class to so-called "nearly-incompressible" vector fields \cite{bianchini2020uniqueness}, which allows the vector field to have unbounded divergence, and a partial extension to vector fields with bounded deformation \cite{ambrosio2005traces}, for which only the symmetric part of the gradient has bounded variation. There is also the work \cite{bianchini2022differentiability}, which gives a notion of differentiability for the flow map associated with a $BV$ vector field under the same assumptions as in \cite{bianchini2020uniqueness}.

    Throughout this paper, we assume that the vector field is divergence-free and time-independent to help illustrate the methods while simplifying the presentation. The complete calculations for a time-dependent vector field can be found in the PhD thesis \cite{huysmans_2024}.

    \setcounter{subsection}{0}
    \subsection{Motivation}
        % There is no analogue of statement \ref{} for the DiPerna-Lions commutator when the vector field lies in $u(x,t)\in L_t^1 BV_x$ as shown in \cite{}, see Appendix \ref{} for more details. However, suppose one aims to give quantitative well-posedness estimates in a weak topology. In that case, one does not necessarily need to control the DiPerna-Lions commutator, but instead construct some other regularisation scheme as shown in the pioneering work of Ambrosio \cite{}. In particular, the regularisation scheme of Ambrosio is local in space, which is in contrast to the regularisation scheme of DiPerna-Lions. This requirement of employing different local smooth approximations to the passive scalar on different sets in space is the key reason behind the tetration nature of the estimates we obtain.

        % Inspired by the following non-quantitative argument, we only need to combine a weak type stability estimate with a compactness argument to get uniform control. One may similarly view our results in the Appendix of \cite{our paper} in this light.
        
        The goal is to give \textit{uniform} estimates for the well-posedness/continuous dependence on initial data for transport along $BV$ vector fields. We first observe that the uniqueness of solutions (due to Ambrosio \cite{ambrosio2004transport}, and re-proven in Corollary \ref{cor: proof of Ambrosio}) can be combined with weak compactness of the unit ball in $BV$ to give such a stability estimate, the proof of which is contained in the appendix:
        \begin{proposition}\label{prop:nonquantitativeBV}
            For any $\epsilon>0$, there exists $\delta>0$ depending only on $\epsilon>0$, $\left\|\rho_0\right\|_{L^\infty}$, $T\left\|\nabla u\right\|_{L^1}$, and the dimension $d\ge2$ such that:
            
            If
            \begin{equation*}
                \left\|\rho_0\right\|_{W^{-1,1}}\le\delta,
            \end{equation*}
            then
            \begin{equation*}
                \left\|\rho_T\right\|_{W^{-1,1}}\le\epsilon.
            \end{equation*}
        \end{proposition}

        Indeed, Theorem \ref{thm:quantitativeBV} is precisely a quantification of this proposition. It is essential in the proof of the above that uniqueness holds when $u(x)\in BV$, rather than the more standard result when $u(x)\in W^{1,1}$ by DiPerna-Lions \cite{diperna1989ordinary}; the same compactness argument would fail in $W^{1,1}$ since the unit ball in $W^{1,1}$ is not weakly compact. It is key, therefore, to exploit well-posedness in $BV$, for which there are only two known proofs \cite{ambrosio2004transport,bianchini2020uniqueness}.
        
        Therefore, our starting point is to quantify the well-posedness proof of Ambrosio \cite{ambrosio2004transport}. In Subsection \ref{sec:quatitativeambrosio}, we quantify Ambrosio's local regularisation scheme \cite{ambrosio2004transport} (see also the lecture notes \cite{ambrosio2014continuity}) to give a quantitative proof of well-posedness for transport by $BV$ vector fields in terms of new harmonic estimates in $BV$. In Subsection \ref{sec:harmonic} we then show how continuous dependence of our bounds on the vector field and passive scalar in the weak-$*$ topologies in $BV$ and $L^\infty$ gives uniform estimates of the kind in Proposition \ref{prop:nonquantitativeBV} for an \textit{open cover} of the vector field and passive scalar in $BV$ and $L^\infty$. Finally, in Subsection \ref{sec:combinatorial}, we apply the pigeonhole principle to give
        a finite sub-cover, and hence a single uniform bound of the form in Proposition \ref{prop:nonquantitativeBV}, which quantifies the dependence of $\delta>0$ on $\epsilon>0$, $\left\|\rho_0\right\|_{L^\infty}$, and $T\left\|\nabla u\right\|_{L^1}$. For the detailed proofs, see Section \ref{section3}.

\section{Outline of the proof}\label{section2}
Let $\rho(x,t)\in L_t^\infty L_x^\infty$ be a solution to \eqref{eq:TE}. We wish to control the weak norm:
\begin{equation}\label{eq:phi_T}
    \left\|\rho_T\right\|_{W^{-1,1}} = \sup_{\left\|\phi_T\right\|_{W^{1,\infty}}\le1}\int_{\mathbb{T}^d}\rho_T(x)\phi_T(x)\;dx.
\end{equation}

Throughout, we fix some test function $\phi_T(x)\in W^{1,\infty}(\mathbb{T}^d;\mathbb{R})$, with $\left\|\phi_T\right\|_{W^{1,\infty}}\le1$. We transport this backwards; let $\phi(x,t)\in L_t^\infty L_x^\infty$ be (any) weak solution to the following \textit{backwards} transport equation:
\begin{equation}\label{eq:phi}
    \begin{aligned}
        & \frac{\partial\phi}{\partial t}(x,t) + \nabla\cdot(u(x)\phi(x,t)) = 0, \\
        & \phi(x,T)=\phi_T(x),
    \end{aligned}
\end{equation}
with the bound
\begin{equation*}
    \left\|\phi\right\|_{L_t^\infty L_x^\infty} \le \left\|\phi_T\right\|_{L^\infty}\le1.
\end{equation*}

\subsection{Quantifying the Ambrosio regularisation}\label{sec:quatitativeambrosio}

% Since the proof involves only a finite number of steps, one can hope for a \textit{finite} number of constraints on the initial datum, in the sense of test functions perhaps, which correspond to smallness in the weak topology, e.g. in negative Sobolev spaces as desired.

In the proof of uniqueness by Ambrosio (see \cite{ambrosio2004transport}, and also the lecture notes \cite{ambrosio2014continuity}), one regularises the solution $\rho(x,t)$, so that one may apply the classical techniques for smooth solutions. One then deduces, as for smooth solutions, that the zero initial datum can only correspond to the zero solution. Careful quantification of this argument should instead translate smallness constraints on the initial datum to smallness constraints on the final datum.

To this end, we find it is easier to work with an adapted version of the proof of Ambrosio, wherein we show instead that the product $\rho(x,t)\phi(x,t)$ is a solution to the transport equation, and so conserves mass. In particular, if
\begin{equation}\label{eq:initialdatumcondition}
    \int_{\mathbb{T}^d}\rho_0(x)\phi(x,0)\;dx=0,
\end{equation}
then
\begin{equation}\label{eq:finaldatumcondition}
    \int_{\mathbb{T}^d}\rho_T(x)\phi(x,T)\;dx=0,
\end{equation}
which controls exactly the weak norm \eqref{eq:phi_T}. To prove this, we now have a choice. To apply the derivative to the product $\rho(x,t)\phi(x,t)$, we need only one of the functions to be smooth. That is, we have a choice to either regularise the solution $\rho(x,t)$ (as in \cite{ambrosio2004transport}) or the solution $\phi(x,t)$. We choose the latter because, once regularised, the condition on the initial datum \eqref{eq:initialdatumcondition} now becomes an integral of $\rho_0(x)$ against a \textit{smooth} function, which one can more easily quantify by smallness in $W^{-1,1}$.

Referring to Section \ref{section3} for the exact expressions and computations, our quantitative regularisation splits the integral \eqref{eq:finaldatumcondition} into five terms,
\begin{equation*}
    \int_{\mathbb{T}^d} \rho_T(x)\phi(x,T)\;dx = I_1 + I_2 + I_3 + I_4 + I_5,
\end{equation*}
each of which we control explicitly. The regularisation itself contains two key ingredients: mollification by an anisotropic mollifier (the Alberti mollifier in \cite{ambrosio2014continuity}) designed to cancel the transport structure and localisation since the Alberti mollifier depends on space.

To be precise, we fix a standard mollifier $\varphi(h)\in C_c^\infty(\mathbb{R}^d;\mathbb{R})$ satisfying certain technical assumptions (we refer to equation \eqref{eq:mollifierassumptions} in Section \ref{section3}). The following three parameters then define the regularisation scheme:

\begin{itemize}
    \item $\Lambda>0$ controls the anisotropic structure of the Alberti mollifiers in \cite{ambrosio2014continuity}, see \eqref{eq:ambrosiomollifier};
    \item $\delta>0$ controls the radius of the Alberti mollifiers, see \eqref{eq:ambrosiomollifierrescaled};
    \item $\epsilon>0$ is the localisation parameter, namely the radius of the spatial ball on which to apply each Alberti mollifier, see \eqref{eq:firstapproximation}.
\end{itemize}

Once we fix these parameters, the five terms arise as follows:
\begin{itemize}
    \item $I_1$ \eqref{eq:I1} is the error from regularising $\phi(x,t)$ in the integral \eqref{eq:finaldatumcondition}, see \eqref{eq:firstapproximation}, and vanishes as $\delta\to0$;
    \item $I_2$ \eqref{eq:I2} is the regularised integral \eqref{eq:initialdatumcondition} after transporting the regularised product $\rho(x,t)\phi(x,t)$, see \eqref{eq:remainingterms}, and vanishes as $\left\|\rho_0\right\|_{W^{-1,1}}\to0$.
\end{itemize} 
    The remaining terms now come from the transport of $\phi(x,t)\rho(x,t)$, and due to the regularisation, we are left with three errors;
\begin{itemize}
    \item $I_3$ \eqref{eq:I3} is the transport error from the localisation and vanishes as $\delta\to0$ for fixed $\epsilon>0$;
    \item $I_4+I_5$ \eqref{eq:I4I5} is the DiPerna-Lions-style transport error (see also \cite{diperna1989ordinary}) from the mollification of $\phi(x,t)$. We split this final error into two terms. The first, $I_4$ \eqref{eq:I4}, vanishing by localisation ($\epsilon\to0$). The second, $I_5$ \eqref{eq:I5}, vanishes by the structure of the Alberti mollifier \eqref{eq:alberticancellation} ($\Lambda\to\infty$).
\end{itemize}

These five terms then have the following explicit bounds in terms of the solution $\phi(x,t)$ in \eqref{eq:phi}:
\begin{proposition}\label{prop:fiveterms}
    For fixed $\Lambda, \delta, \epsilon>0$ with $\delta\le \epsilon e^{-\Lambda}$:
    \begin{align*}
        & |I_1| \lesssim \delta e^{(d+1)\Lambda}\left\|\rho\right\|_{L_t^\infty L_x^\infty}, \\
        & |I_2| \lesssim \bigg(1 + \frac{e^\Lambda}{\delta} + \frac{1}{\epsilon}\bigg)\left\|\rho_0\right\|_{W^{-1,1}}, \\
        & |I_3| \lesssim T^\frac{d-1}{d} \frac{e^{d\Lambda}}{\epsilon} \left\|\rho\right\|_{L_t^\infty L_x^\infty} \left\|\nabla u\right\|_{L^1}\left(\left\|\phi\right\|_{L_t^2L_x^2}^2 - \left\|\phi*\varphi_{2\delta e^\Lambda}\right\|_{L_t^2L_x^2}^2 \right)^\frac{1}{d}, \\
        & |I_4| \lesssim Te^{2\Lambda}\left\|\rho\right\|_{L_t^\infty L_x^\infty} \left\|\nabla u\right\|_{L^1}^\frac{1}{2}\left(\left\|\nabla u\right\|_{L^1} - \left\|\nabla u * \varphi_{4\epsilon}\right\|_{L^1} \right)^\frac{1}{2}, \\
        & |I_5| \lesssim \frac{T}{\Lambda}\left\|\rho\right\|_{L_t^\infty L_x^\infty}\left\|\nabla u\right\|_{L^1},
    \end{align*}
    where $\lesssim$ denotes less than or equal to up to a positive constant depending only on the dimension $d\ge2$ and the choice of mollifier $\varphi(h)\in C_c^\infty(\mathbb{R}^d;\mathbb{R})$ satisfying \eqref{eq:mollifierassumptions}.
\end{proposition}

\subsection{Harmonic estimates in $BV$}\label{sec:harmonic}
The novelty here lies in the bound on $I_4$. Namely the high-frequency control on $u(x)\in BV$ through the difference of norms
\begin{equation}\label{eq:BVharmonic}
    \left\|\nabla u\right\|_{L^1}-\left\|\nabla u * \varphi_{4\epsilon}\right\|_{L^1},
\end{equation}
which vanishes as $\epsilon \to 0$ for $u(x)\in BV$, see for instance \cite{evans2018measure}. In contrast, the more standard high-frequency control $\left\|\nabla u - \nabla u * \varphi_\epsilon\right\|_{L^1}$ appearing in the DiPerna-Lions well-posedness theory \cite{diperna1989ordinary} vanishes only for $u(x)\in W^{1,1}$. Therefore, the above difference of norms is a much weaker harmonic estimate. It emerges by approximating only the \textit{unit} part of the gradient $\frac{\nabla u(x)}{|\nabla u(x)|}$ in Ambrosio's approximation scheme, see \eqref{eq:I4}, instead of the full gradient $\nabla u(x)$ as is necessary for the DiPerna-Lions theory.

In particular, Proposition \ref{prop:fiveterms} can be viewed as a quantification of the stability estimate Proposition \ref{prop:nonquantitativeBV} for fixed $\rho(x,t)\in L_t^\infty L_x^\infty$ and $u(x)\in BV$; one may take $\Lambda >0$ large, $\epsilon > 0$ small, and $\delta>0$ small (and so $\left\|\rho_0\right\|_{W^{-1,1}}$ small) to make all five terms in Proposition \ref{prop:fiveterms} vanish. As shown in the appendix, this is enough to reprove the well-posedness of transport by $BV$ vector fields, originally due to Ambrosio \cite{ambrosio2004transport}:

\begin{corollary}\label{cor: proof of Ambrosio}
    If $\rho_0(x)=0$ is the zero initial datum, then $\rho(x,t)=0$ is the zero solution.
\end{corollary}

The aim is to quantify Proposition \ref{prop:nonquantitativeBV} \textit{uniformly}; we should choose $\Lambda, \epsilon, \delta>0$ to make all the terms in Proposition \ref{prop:fiveterms} small depending only on the norms $\left\|\rho\right\|_{L_t^\infty L_x^\infty}$, and $T\left\|\nabla u\right\|_{L^1}$. To achieve this, the bounds on $I_3, I_4$ in terms of the harmonic estimates on $\phi(x,t) \in L_t^2L_x^2$ and $\nabla u(x)\in BV$, are insufficient.

Consider for instance the harmonic estimate on $\phi(x,t)$ in the bound on $I_3$:
\begin{equation}\label{eq:harmonicestimate}
    \left\|\phi\right\|_{L_t^2L_x^2}^2 - \left\|\phi*\varphi_{2\delta e^\Lambda}\right\|_{L_t^2L_x^2}^2.
\end{equation}
In terms of the Fourier series for $\phi(\cdot,t):\mathbb{T}^d\to\mathbb{R}$, this harmonic estimate roughly corresponds to the $L_t^2L_\omega^2$-norm in frequency-space ($\omega\in\mathbb{Z}^d$) of the frequencies of $\phi(\cdot,t)$ above the threshold $|\omega|\ge\frac{1}{2\delta e^\Lambda}$. The difficulty here is that one has no a-priori knowledge of how high (in frequency-space) the frequencies of $\phi(\cdot,t)$ are located.

To this end, we revisit the terms $I_3$ \eqref{eq:I3} and $I_4+I_5$ \eqref{eq:I4I5}. We notice that these expressions are continuous in the weak-$*$ topologies of $\phi(x,t)\in L_t^2L_x^2$ and $u(x)\in BV$, and so we look for new bounds that respect this continuity. To do so, we perform an (even higher) frequency cutoff on $\phi(x,t)\in L_t^2L_x^2$ and $u(x)\in BV$ to extract a remainder term which does not depend on these even higher frequencies. Splitting them into two parts, a \textit{finite} frequency part and an `even higher' frequency remainder, we have:
\begin{align*}
    I_3 = I_3' + R_1, \\
    I_4 + I_5 = I_4' + I_5' + R_2,
\end{align*}
where
\begin{itemize}
    \item $I_3'$ is bounded by the \textit{finite} range of frequencies
    \begin{equation*}
        \left\|\phi * \varphi_{\delta'}\right\|_{L_t^2L_x^2}^2 - \left\|\phi * \varphi_{\delta'} * \varphi_{2\delta e^\Lambda}\right\|_{L_t^2L_x^2}^2,
    \end{equation*}
    for a free parameter $\delta'>0$ which defines a second/higher frequency cutoff $(\phi * \varphi_{\delta'})(x,t)$.
    \item $I_4'$ is analogously bounded by the $BV$ harmonic estimate
    \begin{equation*}
        \left\|\nabla u * \varphi_{\epsilon'}\right\|_{L^1} - \left\|\nabla u * \varphi_{\epsilon'} * \varphi_{4\epsilon}\right\|_{L^1},
    \end{equation*}
    for a free parameter $\epsilon'>0$ which defines the `even higher' frequency cutoff $(\nabla u * \varphi_{\epsilon'})(x,t)$.
    \item $I_5'$ maintains the same bound as $I_5$.
    \item The `higher' frequency remainders $R_1, R_2$ are small when $\delta',\epsilon'>0$ are small, independent of the `even higher' frequencies of $\phi(x,t)$ and $\nabla u(x)$:
\end{itemize}

\begin{proposition}\label{prop:seventerms}
    For fixed $\Lambda,\delta,\delta',\epsilon,\epsilon'>0$ with $\delta\le\epsilon e^{-\Lambda}$, and a mollifier $\varphi(h)\in C_c^\infty(\mathbb{R}^d;\mathbb{R})$ satisfying the assumptions \eqref{eq:mollifierassumptions}:
    \begin{equation*}
        \int_{\mathbb{T}^d}\rho_T(x)\phi_T(x)\;dx = I_1 + I_2 + I_3' + I_4' + I_5' + R_1 + R_2,
    \end{equation*}
    with the bounds:
    \begin{align*}
        & |I_1| \lesssim \delta e^{(d+1)\Lambda}\left\|\rho\right\|_{L_t^\infty L_x^\infty}, \\
        & |I_2| \lesssim \bigg(1 + \frac{e^\Lambda}{\delta} + \frac{1}{\epsilon}\bigg)\left\|\rho_0\right\|_{W^{-1,1}}, \\
        & |I_3'| \lesssim T^\frac{d-1}{d} \frac{e^{d\Lambda}}{\epsilon} \left\|\rho\right\|_{L_t^\infty L_x^\infty} \left\|\nabla u\right\|_{L^1}\left(\left\|\phi*\varphi_{\delta'}\right\|_{L_t^2L_x^2}^2 - \left\|\phi*\varphi_{\delta'}*\varphi_{2\delta e^\Lambda}\right\|_{L_t^2L_x^2}^2 \right)^\frac{1}{d}, \\
        & |I_4'| \lesssim Te^{2\Lambda}\left\|\rho\right\|_{L_t^\infty L_x^\infty} \left\|\nabla u\right\|_{L^1}^\frac{1}{2} \left(\left\|\nabla u*\varphi_{\epsilon'}\right\|_{L^1} - \left\|\nabla u * \varphi_{\epsilon'} * \varphi_{4\epsilon}\right\|_{L^1} \right)^\frac{1}{2}, \\
        & |I_5'| \lesssim \frac{T}{\Lambda}\left\|\rho\right\|_{L_t^\infty L_x^\infty}\left\|\nabla u\right\|_{L^1}, \\
        & |R_1| \lesssim T\frac{\delta'e^{(d+1)\Lambda}}{\epsilon\delta} \left\|\rho\right\|_{L_t^\infty L_x^\infty} \left\|\nabla u\right\|_{L^1}, \\
        & |R_2| \lesssim T \frac{\epsilon' e^{(d+1)\Lambda}}{\delta} \left\|\rho\right\|_{L_t^\infty L_x^\infty} \left\|\nabla u\right\|_{L^1}.
    \end{align*}
\end{proposition}

% The problematic terms are still $I_3'$ and $I_4'$, but the extra degree of freedom granted by $\delta', \epsilon'>0$ gives us a quantification of Proposition \ref{prop:nonquantitativeBV} not just for fixed $\phi(x,t)$ and $u(x)$, but for a \textit{weak}-$*$ open neighbourhood of $\phi(x,t)\in L_t^2L_x^2$ and $\nabla u(x)\in L^1$. By the Banach–Alaoglu theorem there is then hope to find a finite cover, i.e. a finite set $\{(\Lambda_n,\delta_n,\delta'_n,\epsilon_n,\epsilon'_n)\}_{n=1}^N$ such that the bounds in Proposition \ref{prop:seventerms} are always small for at least one these choices of parameters. This is our strategy.

The problematic terms are still $I_3'$ and $I_4'$, but by the extra degree of freedom granted by $\delta', \epsilon'>0$, for any choice of $\Lambda,\delta,\delta',\epsilon,\epsilon'>0$ which makes the bounds in Proposition \ref{prop:seventerms} suitably small, the \textit{same} choice of parameters is also sufficient for nearby $(\phi*\varphi_{\delta'})(x,t)\in L_t^2L_x^2$ and $(\nabla u * \varphi_{\epsilon'})(x)\in L^1$, i.e. for open neighbourhoods in the \textit{weak-$*$} topology. By the Banach–Alaoglu theorem there is then hope to find a finite cover, i.e. a finite set $\{(\Lambda_n,\delta_n,\delta'_n,\epsilon_n,\epsilon'_n)\}_{n=1}^N$ such that the bounds in Proposition \ref{prop:seventerms} are always small for at least one these choices of parameters. This is our strategy.

% The problematic terms are still $I_3'$ and $I_4'$, but the extra degree of freedom granted by $\delta', \epsilon'>0$ gives us a quantification of Proposition \ref{prop:nonquantitativeBV} not just for fixed $\phi(x,t)$ and $u(x)$, but for a \textit{weak}-$*$ open neighbourhood of $\phi(x,t)\in L_t^2L_x^2$ and $\nabla u(x)\in L^1$; for each $\phi(x,t)$ and $u(x)$ there exists some choice of parameters $\Lambda,\delta,\delta',\epsilon,\epsilon'>0$ to make the bounds in Proposition \ref{prop:seventerms} suitably small, and the \textit{same} choice of $\Lambda,\delta,\delta',\epsilon,\epsilon'>0$ is also sufficient for nearby $(\phi*\varphi_{\delta'})(x,t)\in L_t^2L_x^2$ and $(\nabla u * \varphi_{\epsilon'})(x)\in L^1$, i.e. for open neighbourhoods in the \textit{weak-$*$} topology. By the Banach–Alaoglu theorem there is then hope to find a finite cover, i.e. a finite set $\{(\Lambda_n,\delta_n,\delta'_n,\epsilon_n,\epsilon'_n)\}_{n=1}^N$ such that the bounds in Proposition \ref{prop:seventerms} are always small for at least one these choices of parameters. This is our strategy.

\subsection{The pigeonhole principle}\label{sec:combinatorial}

Consider the bound on $I_3'$. Keeping only the dependence on $\delta, \delta',\epsilon$, we approximately require that the following harmonic estimate is small:
\begin{equation}\label{eq:harmonicestimate2}
    \left\|\phi*\varphi_{\delta'}\right\|_{L_t^2L_x^2}^2 - \left\|\phi*\varphi_{\delta'}*\varphi_{\delta}\right\|_{L_t^2L_x^2}^2 \ll \epsilon^d.
\end{equation}

In terms of the Fourier series for $\phi(\cdot,t):\mathbb{Z}^d\to\mathbb{R}$, this harmonic estimate roughly corresponds to the $L_t^2L_\omega^2$-norm in frequency-space ($\omega\in\mathbb{Z}^d$) of the frequencies of $\phi(\cdot,t)$ in the frequency band $\frac{1}{\delta_n}\le|\omega|\le \frac{1}{\delta'_n}$. A-priori, as for the harmonic estimate on $I_3$ \eqref{eq:harmonicestimate}, we do not know where the frequencies of $\phi(\cdot,t)$ are located. However, for a finite cover $\{(\Lambda_n,\delta_n,\delta'_n,\epsilon_n,\epsilon'_n)\}_{n=1}^N$ then perhaps the different frequency bands $\left[\frac{1}{\delta_n},\frac{1}{\delta'_n}\right]$ are \textit{disjoint}. By the \textit{pigeonhole principle} there is then at least one \textit{empty} pigeonhole, i.e. one frequency band $\left[\frac{1}{\delta_n},\frac{1}{\delta'_n}\right]$ for which the harmonic estimate \eqref{eq:harmonicestimate2} is small with the smallness determined by the number of pigeonholes and the total norm of $\left\|\phi\right\|_{L_t^2L_x^2}^2$. Referring to Subsection \ref{proof:quantitativeBV} for the exact computations, we now apply the pigeonhole principle in this way to find $\delta,\delta',\epsilon,\epsilon'>0$ for which \textit{both} $I_3'$ and $I_4'$ are small, independent of $\phi(x,t)$ and $\nabla u(x)$ (except their norms $\left\|\phi\right\|_{L_t^2L_x^2}$ and $\left\|\nabla u\right\|_{L^1}$).

Firstly, the parameter $\delta'>0$ must be sufficiently small to make the bound on the remainder $R_1$ small but is otherwise free. Keeping only the dependence on $\delta,\epsilon$, this approximately requires
\begin{align*}
    & \delta' \approx \epsilon \delta,
\end{align*}
and so the harmonic estimate \eqref{eq:harmonicestimate2} on $I_3'$ can be rewritten as
\begin{align*}
    \left\|\phi*\varphi_{\epsilon\delta}\right\|_{L_t^2L_x^2}^2 - \left\|\phi*\varphi_{\epsilon\delta}*\varphi_{\delta}\right\|_{L_t^2L_x^2}^2 & \ll \epsilon^d.
\end{align*}
Fixing some $\epsilon>0$ to be chosen later, disjointness of the frequency bands $\left[\frac{1}{\delta_n},\frac{1}{\epsilon\delta_n}\right]$ can be ensured by repeatedly replacing
\begin{gather*}
    \delta_{n+1}=\epsilon\delta_n,
\end{gather*}
The number of such pigeonholes required is determined by the smallness, namely $\epsilon^d$. If we perform the replacement $\delta\mapsto\epsilon\delta$ up to $\epsilon^{-d}$ times, we are guaranteed to find such an `empty' pigeonhole. That is, for fixed $\epsilon>0$, we can guarantee $I_3', R_1$ are small by taking, at worst,
\begin{equation*}
    \delta \approx \epsilon^{\epsilon^{-d}}.
\end{equation*}

We now turn to the second problematic term $I_4'$, which is controlled by the $BV$ harmonic estimate
\begin{equation*}
    \left\|\nabla u*\varphi_{\epsilon'}\right\|_{L^1} - \left\|\nabla u * \varphi_{\epsilon'} * \varphi_{4\epsilon}\right\|_{L^1}.
\end{equation*}
We treat this similarly. Firstly, the parameter $\epsilon'>0$ must be sufficiently small to make the bound on the remainder $R_2$ vanish but is otherwise free. One approximately requires
\begin{align*}
    \epsilon' & \approx \delta.
\end{align*}
Since (at worst) $\delta \approx \epsilon^{\epsilon^{-d}}$, $I_4'$ is (at worst) controlled by the difference of norms
\begin{equation*}
    \left\|\nabla u * \varphi_{\epsilon^{\epsilon^{-d}}}\right\|_{L^1} - \left\|\nabla u * \varphi_{\epsilon^{\epsilon^{-d}}} * \varphi_{\epsilon}\right\|_{L^1} \ll 1
\end{equation*}

The argument is again identical, with disjointness of the above differences/pigeonholes ensured by repeatedly replacing
\begin{equation}\label{eq:exponentiation}
    \epsilon_{n+1} = \epsilon_n^{\epsilon_n^{-d}}.
\end{equation}
The key observation here is that, while before we repeated the \textit{multiplication} $\delta\mapsto\epsilon\delta$ to produce \textit{exponentiation} $\epsilon\mapsto \epsilon^{\epsilon^{-d}}$, we must now repeat \textit{exponentiation} to produce \textit{tetration} in the final bound, see Subsection \ref{proof:quantitativeBV} for the exact computations:

{\renewcommand{\thesubsection}{\ref{thm:quantitativeBV}}
\begin{theorem}
    Let $\kappa>0$ for constants $A,B>0$ depending only on the dimension $d\ge2$, if
    \begin{equation*}
        \left\|\rho_0\right\|_{W^{-1,1}} \le \kappa \left(\exp^{\left\lceil A\exp\left(\frac{B}{\kappa}\left\|\rho\right\|_{L_t^\infty L_x^\infty}T\left\|\nabla u\right\|_{L^1}\right)\right\rceil}\left(\frac{1}{\kappa}\left\|\rho\right\|_{L_t^\infty L_x^\infty}\right)\right)^{-1},
    \end{equation*}
    where $\exp^n(x)$ refers to repeated exponentiation $\exp(\exp(\dots(\exp(x))))$ (tetration), then
    \begin{equation*}
        \left\|\rho_T\right\|_{W^{-1,1}} \le \kappa.
    \end{equation*}
\end{theorem}
}
% {\renewcommand{\thesubsection}{\ref{thm:mixingBV}}
% \begin{corollary}
%     For constants $A,B>0$ depending only on the dimension $d\ge2$, then
%     \begin{equation*}
%         \left\|\rho_T\right\|_{W^{-1,1}} \ge \left\|\rho_0\right\|_{W^{-1,1}} \left(\exp^{\left\lceil A\exp\left(B\frac{\left\|\rho_0\right\|_{L^\infty}}{\left\|\rho_0\right\|_{W^{-1,1}}}T\left\|\nabla u\right\|_{L^1}\right)\right\rceil}\left(\frac{\left\|\rho_0\right\|_{ L^\infty}}{\left\|\rho_0\right\|_{W^{-1,1}}}\right)\right)^{-1},
%     \end{equation*}
%     where $\exp^n(x)$ refers to repeated exponentiation $\underbrace{\exp(\exp(\dots(\exp(x))))}_{n}$ (tetration).
% \end{corollary}
% }

\begin{remark}\label{remark:doubleexponential}
    The tetration can be attributed to the factor of $\frac{1}{\epsilon}$ in the bound for $I_3'$. This leads to the bound $\epsilon^{-d}$ in the harmonic estimate \eqref{eq:harmonicestimate2}, and therefore to repeated exponentiation $\epsilon \mapsto \epsilon^{\epsilon^{-d}}$ \eqref{eq:exponentiation} when applying the pigeonhole principle. This factor emerges from the \textit{local} in space nature of Ambrosio's regularisation scheme \cite{ambrosio2004transport}. To be precise, the Alberti mollifier $\varphi_\delta^{\Lambda,\bar{M}(\bar{x})}(h) \in C_c^\infty(\mathbb{R}^d;\mathbb{R})$ (see equation \eqref{eq:ambrosiomollifier}) is used to mollify the solution $\phi(x,t)$ so that the transport error
    \begin{equation*}
        r(x,t) =\left(\frac{\partial}{\partial t} + u(x)\cdot\nabla\right)\left(\phi*\varphi_\delta^{\Lambda,\bar{M}(\bar{x})}\right)(x,t)
    \end{equation*}
    is small \textit{locally}, i.e. for some ball $x\in B_\epsilon(\bar{x})\subset\mathbb{T}^d$ with the mollifier then having to depend on $\bar{x}\in\mathbb{T}^d$. To ensure smallness globally, we perform a cutoff via the mollifier $\varphi_\epsilon(h)$ in equation \eqref{eq:firstapproximation}, whose regularity $\left\|\varphi_\epsilon\right\|_{W^{1,1}}$ produces the factor $\frac{1}{\epsilon}$ appearing in the bound on $I_3$ \eqref{eq:localityfactor}, and similarly for $I_3'$.

    We remark that when $u(x)\in W^{1,p}$ the standard regularisation due to DiPerna-Lions instead vanishes \textit{globally} (e.g. in $L_t^1 L_x^1$, \cite{diperna1989ordinary}). When $p>1$, a similar approach to the one presented in this paper then gives a polynomial-exponential mixing estimate instead of tetration, see \cite[Section 4]{huysmans2024quantitative} by the same authors. We argue, therefore, that the locality of the regularisation scheme and the resulting mixing rate are intimately related.

    % The regularisation $\{\rho_\epsilon\}_{\epsilon>0}\subset L_t^\infty L_x^\infty$ is chosen depending on a choice of $\bar{x}\in\mathbb{T}^d$ so that the transport error
    % \begin{equation*}
    %     r_\epsilon(x,t)=\left(\frac{\partial}{\partial t} + u(x)\cdot\nabla\right)\rho_\epsilon(x,t),
    % \end{equation*}
    % vanishes as $\epsilon\to0$ \textit{locally}, i.e. for $x\in\mathbb{T}^d$ near $\bar{x}\in\mathbb{T}^d$. If $r_\epsilon(x,t)$ were to vanish \textit{globally} (e.g. in $L_t^1L_x^1$), then we claim that the same framework established in this paper gives a double exponential stability estimate and corresponding mixing bound instead.

    % If one could instead regularise $\phi(x,t)$ with a globally small transport error, localisation would be unnecessary, and this factor would not appear. Thus, the mixing rate and the locality of the regularisation scheme are fundamentally connected.
\end{remark}

\section{Main results}\label{section3}

\subsection{Proof of Proposition \ref{prop:fiveterms}}\label{proof:fiveterms}

\begin{proof}
    For notational convenience, we assume that $\nabla u(x)\in W^{1,1}$, which is dense in $BV$ with the same norm, and so the same result follows.
        
    Throughout, $\lesssim$ will denote less than or equal to up to a positive constant depending only on the dimension $d\ge2$ and the choice of mollifier $\varphi(h)\in C_c^\infty(\mathbb{R}^d;\mathbb{R})$ satisfying
    \begin{equation}\label{eq:mollifierassumptions}
        \begin{aligned}
            & \varphi(h) \ge 0, \\
            & \varphi(-h) = \varphi(h), \\
            & \varphi(h) = 0 \text{ for } |h|\ge 1, \\
            & \varphi(h) > 0 \text{ for } |h|\le \frac{1}{2}. 
        \end{aligned}
    \end{equation}

    In the statement of the proposition we have fixed some $\phi_T(x)\in W^{1,\infty}$ with
    \begin{equation}\label{eq:compactnormscaling}
        \left\|\phi_T\right\|_{W^{1,\infty}} = \left\|\phi_T\right\|_{L^\infty} + \left\|\nabla\phi_T\right\|_{L^\infty} \le 1,
    \end{equation}
    and $\phi(x,t)\in L_t^\infty L_x^\infty$ is a (unique \cite{ambrosio2004transport}) weak solution to the \textit{backwards} transport equation
    \begin{equation}\label{eq:phitransport}
        \begin{aligned}
            & \frac{\partial\phi}{\partial t}(x,t) + u(x)\cdot\nabla\phi(x,t) = 0, \\
            & \phi(x,T)=\phi_T(x),
        \end{aligned}
    \end{equation}
    with the transport bound
    \begin{equation}
        \left\|\phi\right\|_{L_t^\infty L_x^\infty} \le \left\|\phi_T\right\|_{L^\infty} \le 1, \label{eq:phibound}.
    \end{equation}
    
    Therefore, there exists some weak trace $\phi_0(x)\in L^\infty$, see for instance \cite{ambrosio2014continuity}, with the bound
    \begin{equation}
        \left\|\phi_0\right\|_{L^\infty} \le \left\|\phi\right\|_{L_t^\infty L_x^\infty} \le \left\|\phi_T\right\|_{L^\infty} \le 1. \label{eq:initialphibound}
    \end{equation}
    % so that the following holds:
    
    % For any test function $\chi(x,t)\in C_c^\infty(\mathbb{T}^d\times[0,T];\mathbb{R})$:
    % \begin{equation}\label{eq:testfunctiontransport}
    %     \begin{aligned}
    %         & \int_{\mathbb{T}^d\times[0,T]}\phi(x,t)\left(\frac{\partial\chi}{\partial t}(x,t)+u(x)\cdot\nabla\chi(x,t)\right)\;dxdt \\
    %         & \qquad = \int_{\mathbb{T}^d}\phi_T(x)\chi(x,T)\;dx - \int_{\mathbb{T}^d}\phi_0(x)\chi(x,0)\;dx.
    %     \end{aligned}
    % \end{equation}

    % Ignoring the subtle issue of regularity, if one were to take $\chi(x,t) = \rho(x,t)$ in \eqref{eq:testfunctiontransport}, one sees that
    % \begin{equation}\label{eq:motivation}
    %     \int_{\mathbb{T}^d}\phi_T(x)\rho_T(x)\;dx = \int_{\mathbb{T}^d}\phi_0(x)\rho_0(x)\;dx,
    % \end{equation}
    % in terms of the initial data $\rho_0(x)$, which is exactly the type of control we wish for. However, at present $\phi_0(x)$ has no regularity higher than $L^\infty$ against which we may exploit the weak bound on $\left\|\rho_0\right\|_{W^{-1,1}}$.

    % To this end, we will apply the regularisation scheme of Ambrosio to the transported test function $\phi(x,t)\in L_t^\infty L_x^\infty$, and study the analogue of \eqref{eq:motivation}, with the hope of extracting a weak bound on the RHS due to smoothness of the regularisation.
    
    To apply the local regularisation scheme of Ambrosio, we fix some space-dependent matrix field to be chosen later $\bar{M}(\bar{x}):\mathbb{T}^d\to\mathbb{R}^{d\times d}$ with zero trace $\mathrm{tr}(\bar{M}(\bar{x}))=0$, and bounded $l^2$-norm $|\bar{M}(\bar{x})|\le1$. For some $\Lambda>0$, define Alberti's anisotropic mollifier \cite{ambrosio2014continuity}:
    \begin{equation}\label{eq:ambrosiomollifier}
        \varphi^{\Lambda,\bar{M}(\bar{x})}(h) = \frac{1}{\Lambda}\int_0^\Lambda\varphi\left(\exp(-\lambda \bar{M}^\dag(\bar{x}))\cdot h\right)\;d\lambda \in C_c^\infty(\mathbb{R}^d;\mathbb{R}),
    \end{equation}
    where $\bar{M}^\dag$ is the adjoint $\bar{M}^\dag_{i,j}=\bar{M}_{j,i}$. Observe that    $\varphi^{\Lambda,\bar{M}(\bar{x})}(h)\in L^1$ is a standard mollifier since $\mathrm{tr}(\bar{M})=0$, so $\det (\exp(-\lambda \bar{M}^\dag)) = 1$. We also define the rescaled mollifier by
    \begin{equation}\label{eq:ambrosiomollifierrescaled}
        \varphi_\delta^{\Lambda,\bar{M}(\bar{x})}(h) = \frac{1}{\delta^d}\varphi^{\Lambda,\bar{M}(\bar{x})}\left(\frac{h}{\delta}\right).
    \end{equation}
    
    From the expression \eqref{eq:ambrosiomollifier}, and the preliminary assumptions on the isotropic mollifier \eqref{eq:mollifierassumptions}, we have the following a-priori point-wise bounds:
    \begin{align}
        & \left|\varphi_\delta^{\Lambda,\bar{M}(\bar{x})}(h)\right|\lesssim \delta^{-d}1_{|h|\le\delta e^\Lambda}, \label{eq:mollifierestimate1} \\
        & \left|\nabla\varphi_\delta^{\Lambda,\bar{M}(\bar{x})}(h)\right|\lesssim \delta^{-d-1}e^\Lambda 1_{|h|\le\delta e^\Lambda}.\label{eq:mollifierestimate2}
    \end{align}

    Using that the trace $\mathrm{tr}(\bar{M}(\bar{x}))=0$, and so the determinant $\mathrm{det}(\exp(-\lambda\bar{M}^\dag(\bar{x})))=1$, we then also have the bounds
    \begin{align}
        & \left\|\varphi_\delta^{\Lambda,\bar{M}(\bar{x})}\right\|_{L^1} \lesssim 1, \label{eq:mollifierestimate4} \\
        & \left\|\nabla\varphi_\delta^{\Lambda,\bar{M}(\bar{x})}\right\|_{L^1} \lesssim \frac{e^\Lambda}{\delta}. \label{eq:mollifierestimate5}
    \end{align}
    
    We recall the identity of Alberti \cite{ambrosio2014continuity}:
    \begin{align}
        h\cdot \bar{M}(\bar{x}) \cdot \nabla \varphi^{\Lambda,\bar{M}(\bar{x})}(h) & = \frac{1}{\Lambda} \int_0^\Lambda -\frac{\partial}{\partial \lambda} \varphi\left(\exp(-\lambda \bar{M}^\dag(\bar{x}))\cdot h\right) d\lambda \nonumber \\
        & = \frac{1}{\Lambda}\left(\varphi(h) - \varphi\left(\exp(-\Lambda \bar{M}^\dag(\bar{x}))\cdot h\right)\right), \label{eq:alberticancellation}
    \end{align}
    and so
    \begin{equation}
        \int_{\mathbb{R}^d} \left|h\cdot \bar{M}(\bar{x}) \cdot \nabla \varphi^{\Lambda,\bar{M}(\bar{x})}(h)\right|\;dh\lesssim \frac{1}{\Lambda}. \label{eq:albertidecay}
    \end{equation}
    
    Define also the (DiPerna-Lions) commutator \cite{diperna1989ordinary}:
    \begin{align}
        r\left(u,\phi;\varphi^{\Lambda,\bar{M}(\bar{x})}_\delta\right)(x,t) & = \int_{\mathbb{R}^d} \phi(x-h,t)(u(x)-u(x-h))\cdot\nabla\varphi_\delta^{\Lambda,\bar{M}(\bar{x})}(h)\;dh \label{eq:ambrosiomollifier1} \\
        & = \int_{\mathbb{R}^d} \phi(x-\delta h,t) \left(\int_0^1 h\cdot\nabla u(x-s\delta h)\cdot\nabla\varphi^{\Lambda,\bar{M}(\bar{x})}(h)\;ds\right)\;dh, \label{eq:ambrosiomollifier2}
    \end{align}
    where to integrate $\phi(y):\mathbb{T}^d\to\mathbb{R}$ over $\mathbb{R}^d$, we extend periodically, noting the support of the mollifier to ensure the integral remains finite.

    The key observation of Ambrosio \cite{ambrosio2014continuity} is the similarity between the expressions \eqref{eq:albertidecay} and \eqref{eq:ambrosiomollifier2}, for an appropriate choice of anisotropic matrix field $\bar{M}(\bar{x}):\mathbb{T}^d\to\mathbb{R}^{d\times d}$.

    We perform Ambrosio's local regularisation on balls of radius $\epsilon>0$. We control the locality by the support of a mollifier $\varphi_\epsilon(h)=\epsilon^{-d}\varphi\big(\frac{h}{\epsilon}\big)$:
    \begin{equation}\label{eq:firstapproximation}
        \int_{\mathbb{T}^d}\rho_T(x)\phi_T(x)\;dx = \begin{aligned}[t]
            & \int_{\mathbb{T}^d}\int_{\mathbb{R}^d}\rho_T(x)\left(\phi_T*\varphi_\delta^{\Lambda,\bar{M}(\bar{x})}\right)(x,T)\varphi_\epsilon(\bar{x}-x)\;d\bar{x}\;dx \\
            & + \int_{\mathbb{T}^d}\int_{\mathbb{R}^d}\rho_T(x)\left(\phi_T(x)-\left(\phi_T*\varphi_\delta^{\Lambda,\bar{M}(\bar{x})}\right)(x)\right)\varphi_\epsilon(\bar{x}-x)\;d\bar{x}\;dx,
        \end{aligned}
    \end{equation}
    where we extend $\bar{M}(\bar{x}):\mathbb{T}^d\to\mathbb{R}^{d\times d}$ periodically over $\mathbb{R}^d$ to perform the integral over $\bar{x}\in\mathbb{R}^d$.

    The second of these terms is $I_1$:
    \begin{equation}\label{eq:I1}
        I_1 = \int_{\mathbb{T}^d}\int_{\mathbb{R}^d}\rho_T(x)\left(\phi_T(x)-\left(\phi_T*\varphi_\delta^{\Lambda,\bar{M}(\bar{x})}\right)(x)\right)\varphi_\epsilon(\bar{x}-x)\;d\bar{x}\;dx,
    \end{equation}
    and can be bounded by
    \begin{align*}
        |I_1| & \le \int_{\mathbb{T}^d}\int_{\mathbb{R}^d}|\rho_T(x)|\left|\int_{\mathbb{R}^d}\left(\phi_T(x)-\phi_T(x-h)\right)\varphi_\delta^{\Lambda,\bar{M}(\bar{x})}(h)\;dh\right||\varphi_\epsilon(x-\bar{x})|\;d\bar{x}\;dx \\
        & \lesssim \int_{\mathbb{T}^d}\int_{\mathbb{R}^d}|\rho_T(x)|\left(\int_{\mathbb{R}^d}|h|\left\|\nabla\phi_T\right\|_{L^\infty}\delta^{-d}1_{|h|\le\delta e^\Lambda}\;dh\right)|\varphi_\epsilon(x-\bar{x})|\;d\bar{x}\;dx \\
        & \lesssim \delta e^{(d+1)\Lambda}\left\|\nabla\phi_T\right\|_{L^\infty}\left\|\varphi_\epsilon\right\|_{L^1}\left\|\rho_T\right\|_{L^1} \\
        & \lesssim \delta e^{(d+1)\Lambda}\left\|\rho\right\|_{L_t^\infty L_x^\infty},
    \end{align*}
    using the point-wise bound on $\varphi_\delta^{\Lambda,\bar{M}(\bar{x})}(h)$ \eqref{eq:mollifierestimate1}, the $W^{1,\infty}$ bound on the test function $\phi_T(x)$ \eqref{eq:compactnormscaling}, and the standard $L^1$ bound on the mollifier $\varphi_\epsilon(h)$, as required.

    We now rewrite the remaining term in \eqref{eq:firstapproximation} in terms of the initial data $\rho_0(x)$, $\phi_0(x)$ by transport of $\rho(x,t)$:
    \begin{align}
        & \int_{\mathbb{T}^d}\int_{\mathbb{R}^d}\rho_T(x)\left(\phi_T*\varphi_\delta^{\Lambda,\bar{M}(\bar{x})}\right)(x,T)\varphi_\epsilon(\bar{x}-x)\;d\bar{x}\;dx \nonumber \\
        & \qquad = \begin{aligned}[t]
            & \int_{\mathbb{T}^d\times[0,T]} \rho(x,t)\left(\frac{\partial}{\partial t}+u(x)\cdot\nabla\right)\left(\int_{\mathbb{R}^d}\left(\phi*\varphi_\delta^{\Lambda,\bar{M}(\bar{x})}\right)(x,t)\varphi_\epsilon(\bar{x}-x)\;d\bar{x}\right)\;dxdt \\
            & + \int_{\mathbb{T}^d}\int_{\mathbb{R}^d} \rho_0(x)\left(\phi_0*\varphi_\delta^{\Lambda,\bar{M}(\bar{x})}\right)(x) \varphi_\epsilon(\bar{x}-x)\;d\bar{x}\;dx,
        \end{aligned} \nonumber \\
        & \qquad = \begin{aligned}[t]
            & \int_{\mathbb{T}^d\times[0,T]}\int_{\mathbb{R}^d}\rho(x,t)\left(\frac{\partial\phi}{\partial t}*\varphi_\delta^{\Lambda,\bar{M}(\bar{x})}\right)(x,t) \varphi_\epsilon(\bar{x}-x)\;d\bar{x}\;dxdt \\
            & + \int_{\mathbb{T}^d\times[0,T]}\int_{\mathbb{R}^d} \rho(x,t)u(x)\cdot\nabla\left(\phi*\varphi_\delta^{\Lambda,\bar{M}(\bar{x})}\right)(x,t)\varphi_\epsilon(\bar{x}-x)\;d\bar{x}\;dxdt \\
            & - \int_{\mathbb{T}^d\times[0,T]}\int_{\mathbb{R}^d} \rho(x,t)\left(\phi*\varphi_\delta^{\Lambda,\bar{M}(\bar{x})}\right)(x,t)u(x)\cdot\nabla\varphi_\epsilon(\bar{x}-x)\;d\bar{x}\;dxdt \\
            & + \int_{\mathbb{T}^d}\int_{\mathbb{R}^d} \rho_0(x)\left(\phi_0*\varphi_\delta^{\Lambda,\bar{M}(\bar{x})}\right)(x) \varphi_\epsilon(\bar{x}-x)\;d\bar{x}\;dx,
        \end{aligned} \label{eq:remainingterms}
    \end{align}

    The last of these terms is $I_2$:
    \begin{equation}\label{eq:I2}
        I_2 = \int_{\mathbb{T}^d}\int_{\mathbb{R}^d} \rho_0(x)\left(\phi_0*\varphi_\delta^{\Lambda,\bar{M}(\bar{x})}\right)(x) \varphi_\epsilon(\bar{x}-x)\;d\bar{x}\;dx,
    \end{equation}
    and can be bounded by:
    \begin{align*}
        |I_2| & \le \left\|\rho_0\right\|_{W^{-1,1}}\left\|\int_{\mathbb{R}^d}\left(\phi_0*\varphi_\delta^{\Lambda,\bar{M}(\bar{x})}\right)(x)\varphi_\epsilon(\bar{x}-x)\;d\bar{x}\right\|_{W^{1,\infty}} \\
        & \le \begin{aligned}[t]
            & \left\|\rho_0\right\|_{W^{-1,1}}\left\|\int_{\mathbb{R}^d}\left\|\phi_0*\varphi_\delta^{\Lambda,\bar{M}(\bar{x})}\right\|_{L^\infty}|\varphi_\epsilon(\bar{x}-x)|\;d\bar{x}\right\|_{L^\infty} \\
            & + \left\|\rho_0\right\|_{W^{-1,1}}\left\|\int_{\mathbb{R}^d}\left\|\phi_0*\nabla\varphi_\delta^{\Lambda,\bar{M}(\bar{x})}\right\|_{L^\infty}|\varphi_\epsilon(\bar{x}-x)|\;d\bar{x}\right\|_{L^\infty} \\
            & + \left\|\rho_0\right\|_{W^{-1,1}}\left\|\int_{\mathbb{R}^d}\left\|\phi_0*\varphi_\delta^{\Lambda,\bar{M}(\bar{x})}\right\|_{L^\infty}|\nabla\varphi_\epsilon(\bar{x}-x)|\;d\bar{x}\right\|_{L^\infty}
        \end{aligned} \\
        & \le \begin{aligned}[t]
            & \left\|\rho_0\right\|_{W^{-1,1}}\left\|\phi_0\right\|_{L^\infty}\left\|\varphi_\delta^{\Lambda,\bar{M}(\bar{x})}\right\|_{L^1}\left\|\varphi_\epsilon\right\|_{L^1} \\ & + \left\|\rho_0\right\|_{W^{-1,1}}\left\|\phi_0\right\|_{L^\infty}\left\|\nabla\varphi_\delta^{\Lambda,\bar{M}(\bar{x})}\right\|_{L^1}\left\|\varphi_\epsilon\right\|_{L^1} \\ & + \left\|\rho_0\right\|_{W^{-1,1}}\left\|\phi_0\right\|_{L^\infty}\left\|\varphi_\delta^{\Lambda,\bar{M}(\bar{x})}\right\|_{L^1}\left\|\nabla\varphi_\epsilon\right\|_{L^1}
        \end{aligned} \\
        & \lesssim \left\|\rho_0\right\|_{W^{-1,1}}\left(1+\frac{e^\Lambda}{\delta}+\frac{1}{\epsilon}\right),
    \end{align*}
    where we have used the transport bound on $\left\|\phi_0\right\|_{L^\infty}$ \eqref{eq:initialphibound}, the crude norm bounds on $\varphi_\delta^{\Lambda,\bar{M}(\bar{x})}(h)$ \eqref{eq:mollifierestimate4}, \eqref{eq:mollifierestimate5}, and the standard $L^1$ bounds on the rescaled mollifier $\varphi_\epsilon(h)$, as required.

    We are left to rewrite the remaining terms in \eqref{eq:remainingterms}. To exploit transport of $\phi(x,t)$ in \eqref{eq:phitransport}, we rewrite the first term in \eqref{eq:remainingterms} in terms of the DiPerna-Lions commutator for $\phi(x,t)$ \eqref{eq:ambrosiomollifier1}:
    \begin{align}
        & \int_{\mathbb{T}^d\times[0,T]}\int_{\mathbb{R}^d}\rho(x,t)\left(\frac{\partial\phi}{\partial t}*\varphi_\delta^{\Lambda,\bar{M}(\bar{x})}\right)(x,t) \varphi_\epsilon(\bar{x}-x)\;d\bar{x}\;dxdt \nonumber \\
        & \qquad = - \int_{\mathbb{T}^d\times[0,T]}\int_{\mathbb{R}^d}\rho(x,t)\left(\nabla\cdot(u\phi)*\varphi_\delta^{\Lambda,\bar{M}(\bar{x})}\right)(x,t) \varphi_\epsilon(\bar{x}-x)\;d\bar{x}\;dxdt \nonumber \\
        & \qquad = \begin{aligned}[t]
            & - \int_{\mathbb{T}^d\times[0,T]}\int_{\mathbb{R}^d}\rho(x,t)u(x)\cdot\nabla\left(\phi*\varphi_\delta^{\Lambda,\bar{M}(\bar{x})}\right)(x,t) \varphi_\epsilon(x-\bar{x})\;d\bar{x}\;dxdt \\
            & + \int_{\mathbb{T}^d\times[0,T]}\int_{\mathbb{R}^d}\rho(x,t) r\left(\phi,u;\varphi_\delta^{\Lambda,\bar{M}(\bar{x})}\right)(x,t)\varphi_\epsilon(\bar{x}-x)\;d\bar{x}\;dxdt,
        \end{aligned} \label{eq:firstterm}
    \end{align}
    and we now notice that the first term in \eqref{eq:firstterm} cancels with the second term in \eqref{eq:remainingterms}. Therefore, we can rewrite the remaining terms in \eqref{eq:remainingterms} as:
    \begin{align}
        & \int_{\mathbb{T}^d}\int_{\mathbb{R}^d}\rho_T(x)\phi(x,T)\;dx \nonumber \\
        & \qquad = I_1 + I_2 \nonumber \\
        & \qquad \quad \begin{aligned}[t]
            & - \int_{\mathbb{T}^d\times[0,T]}\int_{\mathbb{R}^d}\rho(x,t)\left(\phi*\varphi_\delta^{\Lambda,\bar{M}(\bar{x})}\right)(x,t) u(x)\cdot\nabla\varphi_\epsilon(\bar{x}-x)\;d\bar{x}\;dxdt \\
            & + \int_{\mathbb{T}^d\times[0,T]}\int_{\mathbb{R}^d}\rho(x,t) r\left(\phi,u;\varphi_\delta^{\Lambda,\bar{M}(\bar{x})}\right)(x,t)\varphi_\epsilon(\bar{x}-x)\;d\bar{x}\;dxdt.
        \end{aligned} \label{eq:remainingterms2}
    \end{align}

    The first of these terms is $I_3$:
    \begin{equation}\label{eq:I3}
        I_3 = - \int_{\mathbb{T}^d\times[0,T]}\int_{\mathbb{R}^d}\rho(x,t)\left(\phi*\varphi_\delta^{\Lambda,\bar{M}(\bar{x})}\right)(x,t) u(x)\cdot\nabla\varphi_\epsilon(\bar{x}-x)\;d\bar{x}\;dxdt,
    \end{equation}
    and for a choice of bounded $\bar{\phi}(x,t):\mathbb{T}^d\times[0,T]\to\mathbb{R}$ to be chosen later, can be rewritten as:
    \begin{equation*}
        I_3 = \begin{aligned}[t]
            & - \int_{\mathbb{T}^d\times[0,T]} \int_{\mathbb{R}^d} \rho(x,t) \left(\left(\phi*\varphi_\delta^{\Lambda,\bar{M}(\bar{x})}\right)(x,t)-\bar{\phi}(x,t)\right) u(x)\cdot\nabla\varphi_\epsilon(\bar{x}-x)\;d\bar{x}\;dxdt \\
            & - \int_{\mathbb{T}^d\times[0,T]} \int_{\mathbb{R}^d} \rho(x,t)\bar{\phi}(x,t) u(x)\cdot\nabla\varphi_\epsilon(\bar{x}-x)\;d\bar{x}\;dxdt,
        \end{aligned}
    \end{equation*}
    where the second term is zero since
    \begin{equation*}
        \int_{\mathbb{R}^d} \nabla\varphi_\epsilon(\bar{x}-x)\;d\bar{x} = 0.
    \end{equation*}

    We may now bound $I_3$ using the point-wise bound on $\varphi_\delta^{\Lambda,\bar{M}(\bar{x})}(h)$ \eqref{eq:mollifierestimate1}:
    \begin{align}
        |I_3| & \le \left\|\rho\right\|_{L_t^\infty L_x^\infty}\int_{\mathbb{T}^d\times[0,T]} \int_{\mathbb{R}^d} \left|\left(\phi*\varphi_\delta^{\Lambda,\bar{M}(\bar{x})}\right)(x,t)-\bar{\phi}(x,t)\right| \left|u(x)\right|\left|\nabla\varphi_\epsilon(\bar{x}-x)\right|\;d\bar{x}\;dxdt \nonumber \\
        & \le \left\|\rho\right\|_{L_t^\infty L_x^\infty}\int_{\mathbb{T}^d\times[0,T]} \int_{\mathbb{R}^d} \left|\int_{\mathbb{R}^d}\left(\phi(x-h,t)-\bar{\phi}(x,t)\right)\varphi_\delta^{\Lambda,\bar{M}(\bar{x})}(h)\;dh\right| \left|u(x)\right|\left|\nabla\varphi_\epsilon(\bar{x}-x)\right|\;d\bar{x}\;dxdt \nonumber \\
        & \lesssim \left\|\rho\right\|_{L_t^\infty L_x^\infty}\int_{\mathbb{T}^d\times[0,T]} \int_{\mathbb{R}^d} \int_{\mathbb{R}^d}\left|\phi(x-h,t)-\bar{\phi}(x,t)\right| \delta^{-d}1_{|h|\le\delta e^\Lambda}\;dh \left|u(x)\right|\left|\nabla\varphi_\epsilon(\bar{x}-x)\right|\;d\bar{x}\;dxdt \nonumber \\
        & \lesssim \frac{1}{\epsilon}\left\|\rho\right\|_{L_t^\infty L_x^\infty}\int_{\mathbb{T}^d\times[0,T]} \int_{\mathbb{R}^d}\left|u(x)\right|\left|\phi(x-h,t)-\bar{\phi}(x,t)\right|\delta^{-d}1_{|h|\le\delta e^\Lambda}\;dh\;dxdt, \label{eq:localityfactor}
    \end{align}
    where we have used the standard $L^1$ bound on the rescaled mollifier $\nabla \varphi_\epsilon(h)$. This bound is still not quite quantitative. By Cauchy-Schwartz:
    \begin{align*}
        |I_3| & \lesssim \begin{aligned}[t]
            & \frac{1}{\epsilon}\left\|\rho\right\|_{L_t^\infty L_x^\infty}\left(\int_{\mathbb{T}^d\times[0,T]} \int_{\mathbb{R}^d} |u(x)|^\frac{d}{d-1} \delta^{-d}1_{|h|\le\delta e^\Lambda}\;dh\;dxdt\right)^\frac{d-1}{d} \\
            & \times \left(\int_{\mathbb{T}^d\times[0,T]}\int_{\mathbb{R}^d}\left|\phi(x-h,t)-\bar{\phi}(x,t)\right|^d \delta^{-d}1_{|h|\le\delta e^\Lambda}\;dh\;dxdt\right)^\frac{1}{d}
        \end{aligned} \\
        & \lesssim \begin{aligned}[t]
            & \frac{1}{\epsilon}\left\|\rho\right\|_{L_t^\infty L_x^\infty}\left(T^\frac{d-1}{d}e^{(d-1)\Lambda}\left\|u\right\|_{L^\frac{d}{d-1}}\right) \\
            & \times \left(\left\|\phi\right\|_{L_t^\infty L_x^\infty}^{d-2}\int_{\mathbb{T}^d\times[0,T]} \int_{\mathbb{R}^d}\left|\phi(x-h,t)-\bar{\phi}(x,t)\right|^2\delta^{-d}1_{|h|\le\delta e^\Lambda}\;dh\;dxdt\right)^\frac{1}{d},
        \end{aligned}
    \end{align*}
    where in the last line, we have also used that $d\ge2$ as assumed in the statement of the proposition.
    
    Since $\varphi(h)\gtrsim 1$ for $|h|\le\frac{1}{2}$ \eqref{eq:mollifierassumptions}:
    \begin{align*}
        \int_{\mathbb{R}^d}\left|\phi(x-h,t)-\bar{\phi}(x,t)\right|^2 \delta^{-d}1_{|h|\le \delta e^\Lambda}\;dh & \lesssim \int_{\mathbb{R}^d}\left|\phi(x-h,t)-\bar{\phi}(x,t)\right|^2 \delta^{-d}(\delta e^\Lambda)^d \varphi_{2\delta e^\Lambda}(h) \;dh \\
        & = e^{d\Lambda}\left((\phi^2*\varphi_{2\delta e^\Lambda})(x,t) -2 (\phi * \varphi_{2\delta e^\Lambda})(x,t)\bar{\phi}(x,t) + \bar{\phi}^2(x,t)\right),
    \end{align*}
    where in the last line, we expand the quadratic expression and integrate over $h\in\mathbb{R}^d$. Taking now $\bar{\phi}(x,t) = (\phi * \varphi_{2\delta e^\Lambda})(x,t)$ we obtain the following bound on $I_3$:
    \begin{align}
        |I_3| & \lesssim T^\frac{d-1}{d}\frac{e^{d\Lambda}}{\epsilon} \left\|\rho\right\|_{L_t^\infty L_x^\infty} \left\|u\right\|_{L^\frac{d}{d-1}}\left\|\phi\right\|_{L_t^\infty L_x^\infty}^\frac{d-2}{d}\left(\int_{\mathbb{T}^d\times[0,T]} (\phi^2*\varphi_{2\delta e^\Lambda})(x,t) - (\phi*\varphi_{2\delta e^\Lambda})^2(x,t)\;dxdt \right)^\frac{1}{d} \nonumber \\
        & \lesssim T^\frac{d-1}{d} \frac{e^{d\Lambda}}{\epsilon} \left\|\rho\right\|_{L_t^\infty L_x^\infty} \left\|\nabla u\right\|_{L^1}\left\|\phi\right\|_{L_t^\infty L_x^\infty}^\frac{d-2}{d}\left(\left\|\phi\right\|_{L_t^2L_x^2}^2 - \left\|\phi*\varphi_{2\delta e^\Lambda}\right\|_{L_t^2L_x^2}^2 \right)^\frac{1}{d} \label{eq:quantitativeharmonic1} \\
        & \lesssim T^\frac{d-1}{d} \frac{e^{d\Lambda}}{\epsilon} \left\|\rho\right\|_{L_t^\infty L_x^\infty} \left\|\nabla u\right\|_{L^1}\left(\left\|\phi\right\|_{L_t^2L_x^2}^2 - \left\|\phi*\varphi_{2\delta e^\Lambda}\right\|_{L_t^2L_x^2}^2 \right)^\frac{1}{d}, \nonumber
    \end{align}
    where we have used the Sobolev embedding on $\nabla u \in L^1$, and in the last line that  $\left\|\phi\right\|_{L_t^\infty L_x^\infty} \le 1$ \eqref{eq:phibound}. Note that the bracket is always positive since $\varphi_{2\delta e^\Lambda}(h)\ge0$ \eqref{eq:mollifierassumptions}.

    We are left to control the last term in \eqref{eq:remainingterms2}. To exploit the relationship between the DiPerna-Lions commutator \eqref{eq:ambrosiomollifier2} and Alberti's anisotropic bound \eqref{eq:albertidecay}, we split this into two final terms:
    \begin{align}
        & \int_{\mathbb{T}^d\times[0,T]}\int_{\mathbb{R}^d}\rho(x,t) r\left(\phi,u;\varphi_\delta^{\Lambda,\bar{M}(\bar{x})}\right)(x,t)\varphi_\epsilon(\bar{x}-x)\;d\bar{x}\;dxdt \label{eq:I4I5} \\
        & \qquad = \int_{\mathbb{T}^d\times[0,T]}\int_{\mathbb{R}^d} \int_{\mathbb{R}^d} \rho(x,t) \phi(x-\delta h,t) \left(\int_0^1 h\cdot\nabla u(x-s\delta h)\cdot\nabla\varphi^{\Lambda,\bar{M}(\bar{x})}(h)\;ds\right)\varphi_\epsilon(\bar{x}-x)\;dh\;d\bar{x}\;dxdt \nonumber \\
        & \qquad = I_4 +I_5, \nonumber
    \end{align}
    where, in terms of the choice of traceless, unit-bounded matrix field $\bar{M}(\bar{x})$:
    \begin{equation}\label{eq:I4}
        I_4 = \begin{aligned}[t]
            & \int_{\mathbb{T}^d\times[0,T]}\int_{\mathbb{R}^d} \int_{\mathbb{R}^d} \rho(x,t) \phi(x-\delta h,t) \left(\int_0^1 h\cdot\left(\nabla u(x-s\delta h) - |\nabla u(x-s\delta h)|\bar{M}(\bar{x})\right)\cdot\nabla\varphi^{\Lambda,\bar{M}(\bar{x})}(h)\;ds\right) \\
            & \qquad \times \varphi_\epsilon(\bar{x}-x)\;dh\;d\bar{x}\;dxdt,
        \end{aligned}
    \end{equation}
    and
    \begin{equation}\label{eq:I5}
        I_5 = \int_{\mathbb{T}^d\times[0,T]}\int_{\mathbb{R}^d} \int_{\mathbb{R}^d} \rho(x,t) \phi(x-\delta h,t) \left(\int_0^1 |\nabla u(x-s\delta h)|h\cdot\bar{M}(\bar{x})\cdot\nabla\varphi^{\Lambda,\bar{M}(\bar{x})}(h)\;ds\right)\varphi_\epsilon(\bar{x}-x)\;dh\;d\bar{x}\;dxdt.
    \end{equation}

    We first bound $I_4$ using the bound on the support of $\nabla\varphi^{\Lambda,\bar{M}(\bar{x})}(h)$ \eqref{eq:mollifierestimate2}, and the standard pointwise bound on the rescaled isotropic mollifier $|\varphi_\epsilon(h)|\lesssim\epsilon^{-d}1_{|h|\le\epsilon}$:
    \begin{align*}
        |I_4| & \le \begin{aligned}[t]
            &\left\|\rho\right\|_{L_t^\infty L_x^\infty}\left\|\phi\right\|_{L_t^\infty L_x^\infty} \int_{\mathbb{T}^d\times[0,T]}\int_{\mathbb{R}^d} \int_{\mathbb{R}^d} \left(\int_0^1 |h|\left|\nabla u(x-s\delta h) - |\nabla u(x-s\delta h)|\bar{M}(\bar{x})\right|\left|\nabla\varphi^{\Lambda,\bar{M}(\bar{x})}(h)\right|\;ds\right) \\
            & \qquad \times \left|\varphi_\epsilon(\bar{x}-x)\right|\;dh\;d\bar{x}\;dxdt
        \end{aligned} \\
        & \lesssim \begin{aligned}[t]
            &\left\|\rho\right\|_{L_t^\infty L_x^\infty} \int_{\mathbb{T}^d\times[0,T]}\int_{\mathbb{R}^d} \int_{\mathbb{R}^d} \int_0^1 |h|\left|\nabla u(x-s\delta h) - |\nabla u(x-s\delta h)|\bar{M}(\bar{x})\right|1_{|h|\le e^\Lambda}\epsilon^{-d}1_{|\bar{x}-x|\le\epsilon} \\
            & \qquad \times \left|\nabla\varphi^{\Lambda,\bar{M}(\bar{x})}(h)\right|\;ds \;dh\;d\bar{x}\;dxdt
        \end{aligned} \\
        & \lesssim \begin{aligned}[t]
            & e^{\Lambda}\left\|\rho\right\|_{L_t^\infty L_x^\infty} \int_{\mathbb{R}^d} \int_0^1\int_{\mathbb{T}^d\times[0,T]}\int_{\mathbb{R}^d} \left|\nabla u(x-s\delta h) - |\nabla u(x-s\delta h)|\bar{M}(\bar{x})\right|\epsilon^{-d}1_{|\bar{x}-x|\le\epsilon}1_{|h|\le e^\Lambda} \\
            & \qquad \times \left|\nabla\varphi^{\Lambda,\bar{M}(\bar{x})}(h)\right|\;d\bar{x}\;dxdt\;ds \;dh
        \end{aligned} \\
        & \lesssim e^{\Lambda}\left\|\rho\right\|_{L_t^\infty L_x^\infty} \int_{\mathbb{R}^d} \int_0^1\int_{\mathbb{T}^d\times[0,T]}\int_{\mathbb{R}^d} \left|\nabla u(x) - |\nabla u(x)|\bar{M}(\bar{x})\right|\epsilon^{-d}1_{|\bar{x}-x|\le\epsilon+\delta e^\Lambda} \left|\nabla\varphi^{\Lambda,\bar{M}(\bar{x})}(h)\right|\;d\bar{x}\;dxdt\;ds \;dh \\
        & \lesssim Te^{\Lambda}\left\|\rho\right\|_{L_t^\infty L_x^\infty} \int_{\mathbb{T}^d}\int_{\mathbb{R}^d} \left|\nabla u(x) - |\nabla u(x)|\bar{M}(\bar{x})\right|\epsilon^{-d}1_{|\bar{x}-x|\le2\epsilon}\left\|\nabla\varphi^{\Lambda,\bar{M}(\bar{x})}\right\|_{L^1}\;d\bar{x}\;dx \\
        & \lesssim Te^{2\Lambda}\left\|\rho\right\|_{L_t^\infty L_x^\infty} \int_{\mathbb{T}^d}\int_{\mathbb{R}^d}|\nabla u(x)|\left|\frac{\nabla u(x)}{|\nabla u(x)|} - \bar{M}(\bar{x})\right|\epsilon^{-d}1_{|\bar{x}-x|\le2\epsilon}\;d\bar{x}\;dx,
    \end{align*}
    where we have used the transport bound on $\phi(x,t)$ \eqref{eq:phibound}, the bound on $\nabla\varphi^{\Lambda,\bar{M}(\bar{x})}(h)$ \eqref{eq:mollifierestimate5}, and the assumption in the statement of the proposition that $\delta e^\Lambda \le \epsilon$. We still face the same problem as for $I_3$, that this bound is still not quantitative. Therefore, by Cauchy-Schwartz:
    \begin{align*}
        |I_4| & \lesssim \begin{aligned}[t]
            & Te^{2\Lambda}\left\|\rho\right\|_{L_t^\infty L_x^\infty} \left(\int_{\mathbb{T}^d} \int_{\mathbb{R}^d} |\nabla u(x)| \epsilon^{-d}1_{|\bar{x}-x|\le2\epsilon}\;d\bar{x}\;dx\right)^\frac{1}{2} \\
            & \times \left(\int_{\mathbb{T}^d}\int_{\mathbb{R}^d} |\nabla u(x)|\left|\frac{\nabla u(x)}{|\nabla u(x)|} - \bar{M}(\bar{x})\right|^2\epsilon^{-d}1_{|\bar{x}-x|\le2\epsilon}\;d\bar{x}\;dx\right)^\frac{1}{2}
        \end{aligned} \\
        & \lesssim Te^{2\Lambda}\left\|\rho\right\|_{L_t^\infty L_x^\infty} \left\|\nabla u\right\|_{L^1}^\frac{1}{2} \left(\int_{\mathbb{T}^d}\int_{\mathbb{R}^d} |\nabla u(x)|\left|\frac{\nabla u(x)}{|\nabla u(x)|} - \bar{M}(\bar{x})\right|^2\epsilon^{-d}1_{|\bar{x}-x|\le2\epsilon}\;dx\;d\bar{x}\right)^\frac{1}{2},
    \end{align*}
    where in the last line, we have also flipped the integrals over $x\in\mathbb{T}^d$, $\bar{x}\in\mathbb{R}^d$ to instead $x\in\mathbb{R}^d$ (extended periodically) and $\bar{x}\in\mathbb{T}^d$. This is valid since $1_{|\bar{x}-x|\le 2\epsilon}$ is symmetric in $\bar{x}-x \in \mathbb{R}^d$.
    
    Since $\varphi(h)\gtrsim 1$ for $|h|\le\frac{1}{2}$ \eqref{eq:mollifierassumptions}:
    \begin{align*}
        & \int_{\mathbb{R}^d} |\nabla u(x)|\left|\frac{\nabla u(x)}{|\nabla u(x)|} - \bar{M}(\bar{x})\right|^2\epsilon^{-d}1_{|\bar{x}-x|\le2\epsilon}\;dx \\
        & \qquad \lesssim \int_{\mathbb{R}^d}|\nabla u(\bar{x}-h)|\left|\frac{\nabla u(\bar{x}-h)}{|\nabla u(\bar{x}-h)|} - \bar{M}(\bar{x})\right|^2 \varphi_{4\epsilon}(h)\;dh \\
        & \qquad = (|\nabla u| * \varphi_{4\epsilon})(\bar{x}) - 2 (\nabla u * \varphi_{4\epsilon})(\bar{x}) \cdot \bar{M}(\bar{x}) + (|\nabla u| * \varphi_{4\epsilon})(\bar{x})|\bar{M}(\bar{x})|^2,
    \end{align*}
    where in the last line, we expanded the quadratic expression and integrated over $h\in\mathbb{R}^d$. Taking now
    \begin{equation*}
        \bar{M}(\bar{x}) = \frac{(\nabla u * \varphi_{4\epsilon})(\bar{x})}{|(\nabla u * \varphi_{4\epsilon})(\bar{x})|},
    \end{equation*}
    which is traceless since $u(x)$ is divergence-free, and by definition unit-bounded, we obtain the following bound on $I_4$:
    \begin{align}
        |I_4| & \lesssim Te^{2\Lambda}\left\|\rho\right\|_{L_t^\infty L_x^\infty} \left\|\nabla u\right\|_{L^1}^\frac{1}{2} \left(\int_{\mathbb{T}^d} 2(|\nabla u| * \varphi_{4\epsilon})(\bar{x}) - 2 |(\nabla u * \varphi_{4\epsilon})(\bar{x})| \;d\bar{x}\right)^\frac{1}{2} \nonumber \\
        & = Te^{2\Lambda}\left\|\rho\right\|_{L_t^\infty L_x^\infty} \left\|\nabla u\right\|_{L^1}^\frac{1}{2} \left(2\left\|\nabla u\right\|_{L^1} - 2\left\|\nabla u * \varphi_{4\epsilon}\right\|_{L^1} \right)^\frac{1}{2}. \label{eq:quantitativeharmonic2}
    \end{align}

    Finally, we exploit Alberti's anisotropic bound \eqref{eq:albertidecay} to bound the remaining term $I_5$ \eqref{eq:I5}:
    \begin{align}
        |I_5| & \le \begin{aligned}[t]
            & \left\|\rho\right\|_{L_t^\infty L_x^\infty}\left\|\phi\right\|_{L_t^\infty L_x^\infty}\int_{\mathbb{T}^d\times[0,T]}\int_{\mathbb{R}^d} \int_{\mathbb{R}^d} \int_0^1 |\nabla u(x-s\delta h)|\left|h\cdot\bar{M}(\bar{x})\cdot\nabla\varphi^{\Lambda,\bar{M}(\bar{x})}(h)\right|\left|\varphi_\epsilon(\bar{x}-x)\right|\; \\
            & \qquad ds\;dh\;d\bar{x}\;dxdt
        \end{aligned} \nonumber \\
        & \lesssim \begin{aligned}[t]
            & \left\|\rho\right\|_{L_t^\infty L_x^\infty} \int_{\mathbb{R}^d} \int_0^1\int_{\mathbb{T}^d\times[0,T]}\int_{\mathbb{R}^d} |\nabla u(x-s\delta h)|\left|h\cdot\bar{M}(\bar{x})\cdot\nabla\varphi^{\Lambda,\bar{M}(\bar{x})}(h)\right|1_{|h|\le e^\Lambda}\epsilon^{-d}1_{|\bar{x}-x|\le\epsilon}\; \\
            & \qquad d\bar{x}\;dxdt\;ds\;dh
        \end{aligned} \nonumber \\
        & \lesssim \begin{aligned}[t]
            & \left\|\rho\right\|_{L_t^\infty L_x^\infty} \int_{\mathbb{R}^d} \int_0^1\int_{\mathbb{T}^d\times[0,T]}\int_{\mathbb{R}^d} |\nabla u(x)|\left|h\cdot\bar{M}(\bar{x})\cdot\nabla\varphi^{\Lambda,\bar{M}(\bar{x})}(h)\right|\epsilon^{-d}1_{|\bar{x}-x|\le\epsilon+\delta e^\Lambda}\; \\
            & \qquad d\bar{x}\;dxdt\;ds\;dh
        \end{aligned} \nonumber \\
        & \lesssim T\left\|\rho\right\|_{L_t^\infty L_x^\infty}  \int_{\mathbb{T}^d} |\nabla u(x)|\epsilon^{-d}\int_{\mathbb{R}^d}1_{|\bar{x}-x|\le\epsilon+\delta e^\Lambda}\int_{\mathbb{R}^d}\left|h\cdot\bar{M}(\bar{x})\cdot\nabla\varphi^{\Lambda,\bar{M}(\bar{x})}(h)\right|\;dh\;d\bar{x}\;dx \nonumber \\
        & \lesssim \frac{T}{\Lambda}\left(1+\frac{\delta e^\Lambda}{\epsilon}\right)^d\left\|\rho\right\|_{L_t^\infty L_x^\infty}\left\|\nabla u\right\|_{L^1} \nonumber \\
        & \lesssim \frac{T}{\Lambda}\left\|\rho\right\|_{L_t^\infty L_x^\infty}\left\|\nabla u\right\|_{L^1}, \label{eq:I5bound}
    \end{align}
    where we have again used the bound on the support of $\nabla\varphi^{\Lambda,\bar{M}(\bar{x})}(h)$ \eqref{eq:mollifierestimate2}, the transport bound on $\phi(x,t)$ \eqref{eq:phibound}, and in the final line the assumption in the statement of the Proposition that $\delta \le \epsilon e^{-\Lambda}$.
\end{proof}

\subsection{Proof of Proposition \ref{prop:seventerms}}\label{proof:seventerms}

\begin{proof}
    We recall the definition of $I_3$, \eqref{eq:I3}:
    \begin{align*}
        I_3 & = - \int_{\mathbb{T}^d\times[0,T]}\int_{\mathbb{R}^d}\rho(x,t)\left(\phi*\varphi_\delta^{\Lambda,\bar{M}(\bar{x})}\right)(x,t) u(x)\cdot\nabla\varphi_\epsilon(\bar{x}-x)\;d\bar{x}\;dxdt \\
        & = I_3' + R_1.
    \end{align*}
    Introducing the second frequency cutoff $\delta'>0$, we split this as follows:
    \begin{equation}\label{eq:I3,1}
        I_3' = - \int_{\mathbb{T}^d\times[0,T]}\int_{\mathbb{R}^d}\rho(x,t)\left(\left(\phi*\varphi_{\delta'}\right)*\varphi_\delta^{\Lambda,\bar{M}(\bar{x})}\right)(x,t) u(x)\cdot\nabla\varphi_\epsilon(\bar{x}-x)\;d\bar{x}\;dxdt,
    \end{equation}
    and
    \begin{equation*}
        R_1 = - \int_{\mathbb{T}^d\times[0,T]}\int_{\mathbb{R}^d}\rho(x,t)\left(\left(\phi-\phi*\varphi_{\delta'}\right)*\varphi_\delta^{\Lambda,\bar{M}(\bar{x})}\right)(x,t) u(x)\cdot\nabla\varphi_\epsilon(\bar{x}-x)\;d\bar{x}\;dxdt.
    \end{equation*}
    $I_3'$ now is the same as $I_3$, with $\phi(x,t)$ replaced with $(\phi*\varphi_{\delta'})(x,t)$, and so is bounded as in \eqref{eq:quantitativeharmonic1} by replacing $\phi(x,t)$ with $(\phi*\varphi_{\delta'})(x,t)$:
    \begin{align*}
        |I_3'| & \lesssim T^\frac{d-1}{d} \frac{e^{d\Lambda}}{\epsilon} \left\|\rho\right\|_{L_t^\infty L_x^\infty} \left\|\nabla u\right\|_{L^1}\left\|\phi*\varphi_{\delta'}\right\|_{L_t^\infty L_x^\infty}^\frac{d-2}{d}\left(\left\|\phi*\varphi_{\delta'}\right\|_{L_t^2L_x^2}^2 - \left\|\phi*\varphi_{\delta'}*\varphi_{2\delta e^\Lambda}\right\|_{L_t^2L_x^2}^2 \right)^\frac{1}{d} \\
        & \lesssim T^\frac{d-1}{d} \frac{e^{d\Lambda}}{\epsilon} \left\|\rho\right\|_{L_t^\infty L_x^\infty} \left\|\nabla u\right\|_{L^1}\left(\left\|\phi*\varphi_{\delta'}\right\|_{L_t^2L_x^2}^2 - \left\|\phi*\varphi_{\delta'}*\varphi_{2\delta e^\Lambda}\right\|_{L_t^2L_x^2}^2 \right)^\frac{1}{d},
    \end{align*}
    where in the last line we have used that $(\phi*\varphi_{\delta'})(x,t)$ inherits the bound $\left\|\phi*\varphi_{\delta'}\right\|_{L_t^\infty L_x^\infty}\le1$ from $\left\|\phi\right\|_{L_t^\infty L_x^\infty}\le1$ \eqref{eq:phibound}.

    We recall now the definition of $I_4 + I_5$, \eqref{eq:I4I5}:
    \begin{align*}
        (I_4 + I_5) & = \int_{\mathbb{T}^d\times[0,T]}\int_{\mathbb{R}^d}\rho(x,t) r\left(\phi,u;\varphi_\delta^{\Lambda,\bar{M}(\bar{x})}\right)(x,t)\varphi_\epsilon(\bar{x}-x)\;d\bar{x}\;dxdt \\
        & = (I_4' + I_5') + R_2.
    \end{align*}
    Introducing the second frequency cutoff $\epsilon'>0$, we split this as follows:
    \begin{equation*}
        (I_4' + I_5') = \int_{\mathbb{T}^d\times[0,T]}\int_{\mathbb{R}^d}\rho(x,t) r\left(\phi,u*\varphi_{\epsilon'};\varphi_\delta^{\Lambda,\bar{M}(\bar{x})}\right)(x,t)\varphi_\epsilon(\bar{x}-x)\;d\bar{x}\;dxdt,
    \end{equation*}
    and
    \begin{equation*}
        R_2 = \int_{\mathbb{T}^d\times[0,T]}\int_{\mathbb{R}^d}\rho(x,t) r\left(\phi,u-u*\varphi_{\epsilon'};\varphi_\delta^{\Lambda,\bar{M}(\bar{x})}\right)(x,t)\varphi_\epsilon(\bar{x}-x)\;d\bar{x}\;dxdt,
    \end{equation*}
    noting that the commutator is linear in the argument $u(x)$ \eqref{eq:ambrosiomollifier1}. $I_4'$ and $I_5'$ are now the same as $I_4$ \eqref{eq:I4} and $I_5$ \eqref{eq:I5}, with $u(x)$ replaced with $(u*\varphi_{\epsilon'})(x)$, and so are bounded as in \eqref{eq:quantitativeharmonic2} and \eqref{eq:I5bound} by replacing $u(x)$ with $(u*\varphi_{\epsilon'})(x)$:
    \begin{align*}
        |I_4'| & \lesssim Te^{2\Lambda}\left\|\rho\right\|_{L_t^\infty L_x^\infty} \left\|\nabla u * \varphi_{\epsilon'}\right\|_{L^1}^\frac{1}{2} \left(\left\|\nabla u*\varphi_{\epsilon'}\right\|_{L^1} - \left\|\nabla u * \varphi_{\epsilon'} * \varphi_{4\epsilon}\right\|_{L^1} \right)^\frac{1}{2} \\
        & \lesssim Te^{2\Lambda}\left\|\rho\right\|_{L_t^\infty L_x^\infty} \left\|\nabla u\right\|_{L^1}^\frac{1}{2} \left(\left\|\nabla u*\varphi_{\epsilon'}\right\|_{L^1} - \left\|\nabla u * \varphi_{\epsilon'} * \varphi_{4\epsilon}\right\|_{L^1} \right)^\frac{1}{2},
    \end{align*}
    and
    \begin{align*}
        |I_5'| & \lesssim \frac{T}{\Lambda}\left\|\rho\right\|_{L_t^\infty L_x^\infty}\left\|\nabla u*\varphi_{\epsilon'}\right\|_{L^1} \\
        & \lesssim \frac{T}{\Lambda}\left\|\rho\right\|_{L_t^\infty L_x^\infty}\left\|\nabla u\right\|_{L^1}.
    \end{align*}

    We now turn to the remainders $R_1$ and $R_2$. First:
    \begin{align*}
        |R_1| & = \left| \int_{\mathbb{T}^d\times[0,T]}\int_{\mathbb{R}^d}\rho(x,t)\left(\left(\phi*\varphi_\delta^{\Lambda,\bar{M}(\bar{x})}\right)(x,t)-\left(\phi*\varphi_\delta^{\Lambda,\bar{M}(\bar{x})}*\varphi_{\delta'}\right)(x,t)\right) u(x)\cdot\nabla\varphi_\epsilon(\bar{x}-x)\;d\bar{x}\;dxdt \right| \\
        & \lesssim \begin{aligned}[t]
            & \int_{\mathbb{T}^d\times[0,T]}\int_{\mathbb{R}^d}|\rho(x,t)|\left|\int_{\mathbb{R}^d}\left(\left(\phi*\varphi_\delta^{\Lambda,\bar{M}(\bar{x})}\right)(x,t)-\left(\phi*\varphi_\delta^{\Lambda,\bar{M}(\bar{x})}\right)(x-h,t)\right)\varphi_{\delta'}(h)\;dh\right| \\
            & \qquad \times |u(x)||\nabla\varphi_\epsilon(\bar{x}-x)|\;d\bar{x}\;dxdt
        \end{aligned} \\
        & \lesssim \begin{aligned}[t]
            & \int_{\mathbb{T}^d\times[0,T]}\int_{\mathbb{R}^d}|\rho(x,t)|\left|\int_{\mathbb{R}^d}\left(\int_0^1 h\cdot\nabla\left(\phi*\varphi_\delta^{\Lambda,\bar{M}(\bar{x})}\right)(x-sh,t) \;ds\right)\varphi_{\delta'}(h)\;dh\right| \\
            & \qquad \times |u(x)||\nabla\varphi_\epsilon(\bar{x}-x)|\;d\bar{x}\;dxdt
        \end{aligned} \\
        & \lesssim \begin{aligned}[t]
            & \int_{\mathbb{T}^d\times[0,T]}\int_{\mathbb{R}^d}\int_{\mathbb{R}^d}\int_0^1|\rho(x,t)|\delta'\left|\left(\phi*\nabla\varphi_\delta^{\Lambda,\bar{M}(\bar{x})}\right)(x-sh,t)\right||\varphi_{\delta'}(h)||u(x)||\nabla\varphi_\epsilon(\bar{x}-x)| \\
            & \qquad \;ds\;dh\;d\bar{x}\;dxdt
        \end{aligned} \\
        & \lesssim \delta' \begin{aligned}[t]
            & \int_{\mathbb{T}^d\times[0,T]}\int_{\mathbb{R}^d}\int_{\mathbb{R}^d}\int_0^1|\rho(x,t)|\left|\int_{\mathbb{R}^d}\phi(x-sh-h',t)\nabla\varphi_\delta^{\Lambda,\bar{M}(\bar{x})}(h')\;dh'\right||\varphi_{\delta'}(h)||u(x)||\nabla\varphi_\epsilon(\bar{x}-x)| \\
            & \qquad \;ds\;dh\;d\bar{x}\;dxdt
        \end{aligned} \\
        & \lesssim \begin{aligned}[t]
            & \delta' \int_{\mathbb{T}^d\times[0,T]}\int_{\mathbb{R}^d}\int_{\mathbb{R}^d}\int_0^1\int_{\mathbb{R}^d}|\rho(x,t)||\phi(x-sh-h',t)|\delta^{-d-1}e^\Lambda 1_{|h'|\le \delta e^\Lambda}|\varphi_{\delta'}(h)||u(x)||\nabla\varphi_\epsilon(\bar{x}-x)| \\
            & \qquad \;dh'\;ds\;dh\;d\bar{x}\;dxdt,
        \end{aligned}
    \end{align*}
    where in the last line, we have used the point-wise bound on the mollifier $\nabla\varphi_\delta^{\Lambda,\bar{M}(\bar{x})}(h)$ \eqref{eq:mollifierestimate2}. Integrating first over $\bar{x}\in\mathbb{R}^d$ and then over $(x,t)\in\mathbb{T}^d\times[0,T]$:
    \begin{align*}
        |R_1| & \lesssim \delta' \left\|\nabla \varphi_\epsilon\right\|_{L^1} \left\|\rho\right\|_{L_t^\infty L_x^\infty} \left\|\phi\right\|_{L_t^1 L_x^d}\left\|u\right\|_{L^\frac{d}{d-1}} \int_{\mathbb{R}^d}\int_0^1\int_{\mathbb{R}^d}\delta^{-d-1}e^\Lambda 1_{|h'|\le \delta e^\Lambda}|\varphi_{\delta'}(h)|\;dh'\;ds\;dh \\
        & \lesssim \delta' \left\|\nabla \varphi_\epsilon\right\|_{L^1} \left\|\rho\right\|_{L_t^\infty L_x^\infty} \left\|\phi\right\|_{L_t^1 L_x^d}\left\|u\right\|_{L^\frac{d}{d-1}}\left\|\varphi_{\delta'}\right\|_{L^1}\delta^{-d-1}e^\Lambda\int_{\mathbb{R}^d}1_{|h'|\le\delta e^\Lambda}\;dh' \\
        & \lesssim \frac{\delta'}{\epsilon} \left\|\rho\right\|_{L_t^\infty L_x^\infty} T \left\|\nabla u\right\|_{L^1}\delta^{-1}e^{(d+1)\Lambda},
    \end{align*}
    where in the last line, we have used the standard bound $\left\|\nabla \varphi_\epsilon\right\|_{L^1}\lesssim \frac{1}{\epsilon}$, that on the torus $\left\|\phi\right\|_{L_t^1 L_x^d}\lesssim T\left\|\phi\right\|_{L_t^\infty L_x^\infty}\le 1$ \eqref{eq:phibound}, the Sobolev embedding $\left\|u\right\|_{L^\frac{d}{d-1}}\lesssim\left\|\nabla u\right\|_{L^1}$, and the standard bound $\left\|\varphi_{\delta'}\right\|_{L^1}\lesssim 1$.

    Finally, for the remainder $R_2$, we first use the expression for the commutator \eqref{eq:ambrosiomollifier1}:
    \begin{align*}
        |R_2| & = \left|\int_{\mathbb{T}^d\times[0,T]}\int_{\mathbb{R}^d}\rho(x,t) r\left(\phi,u-u*\varphi_{\epsilon'};\varphi_\delta^{\Lambda,\bar{M}(\bar{x})}\right)(x,t)\varphi_\epsilon(\bar{x}-x)\;d\bar{x}\;dxdt\right| \\
        & \lesssim \begin{aligned}[t]
            & \int_{\mathbb{T}^d\times[0,T]}\int_{\mathbb{R}^d}|\rho(x,t)|\left|\int_{\mathbb{R}^d} \phi(x-h,t)\left((u-u*\varphi_{\epsilon'})(x)-(u-u*\varphi_{\epsilon'})(x-h)\right)\cdot\nabla\varphi_\delta^{\Lambda,\bar{M}(\bar{x})}(h)\;dh\right| \\
            & \qquad \times |\varphi_\epsilon(\bar{x}-x)|\;d\bar{x}\;dxdt
        \end{aligned} \\
        & \lesssim \begin{aligned}[t]
            & \int_{\mathbb{T}^d\times[0,T]}\int_{\mathbb{R}^d}\int_{\mathbb{R}^d}|\rho(x,t)||\phi(x-h,t)|\left(|u-u*\varphi_{\epsilon'})(x)|+|(u-u*\varphi_{\epsilon'})(x-h)|\right)\delta^{-d-1}e^\Lambda1_{|h|\le\delta e^\Lambda} \\
            & \qquad \times |\varphi_\epsilon(\bar{x}-x)|\;dh\;d\bar{x}\;dxdt,
        \end{aligned}
    \end{align*}
    where in the last line, we have used the point-wise bound on the mollifier $\nabla\varphi_\delta^{\Lambda,\bar{M}(\bar{x})}(h)$ \eqref{eq:mollifierestimate2}. Integrating first over $\bar{x}\in\mathbb{R}^d$, and then over $(x,t)\in\mathbb{T}^d\times[0,T]$:
    \begin{align*}
        |R_2| & \lesssim \left\|\varphi_\epsilon\right\|_{L^1}\left\|\rho\right\|_{L_t^\infty L_x^\infty} \left\|\phi\right\|_{L_t^1L_x^\infty}\left\|u-u*\varphi_{\epsilon'}\right\|_{L^1}\delta^{-d-1}e^\Lambda\int_{\mathbb{R}^d}1_{|h|\le\delta e^\Lambda}\;dh \\
        & \lesssim \left\|\rho\right\|_{L_t^\infty L_x^\infty} T \epsilon' \left\|\nabla u\right\|_{L^1} \delta^{-1}e^{(d+1)\Lambda},
    \end{align*}
    where in the last line we have used the standard bound $\left\|\varphi_\epsilon\right\|_{L^1}\lesssim 1$, that $\left\|\phi\right\|_{L_t^1L_x^\infty}\le T\left\|\phi\right\|_{L_t^\infty L_x^\infty}\le 1$ \eqref{eq:phibound}, and the standard bound $\left\|u-u*\varphi_{\epsilon'}\right\|_{L^1}\lesssim \epsilon'\left\|\nabla u\right\|_{L^1}$.
\end{proof}

\subsection{Proof of Theorem \ref{thm:quantitativeBV} and Corollary \ref{thm:mixingBV}}\label{proof:quantitativeBV}

\begin{proof}
    Throughout this proof, $\lesssim$ will denote less than or equal to up to a positive constant depending only on the dimension $d\ge2$ and the choice of mollifier $\varphi(h)\in C_c^\infty(\mathbb{R}^d;\mathbb{R})$ satisfying \eqref{eq:mollifierassumptions}.
    
    Fix some $\kappa>0$. We use Proposition \ref{prop:seventerms}, that is for any fixed test function $\phi_T(x)\in W^{1,\infty}$ with
    \begin{equation*}
        \left\|\phi_T\right\|_{W^{1,\infty}}\le1,
    \end{equation*}
    we aim to bound:
    \begin{equation}\label{eq:kappabound}
        \int_{\mathbb{T}^d}\rho_T(x)\phi_T(x)\;dx\lesssim \kappa.
    \end{equation}
    To this end, it is sufficient for us to choose parameters as in Proposition \ref{prop:seventerms} so that we bound each of the following seven terms by $\kappa>0$, namely
    \begin{equation}\label{eq:deltabound}
        \delta\le\epsilon e^{-\Lambda},
    \end{equation}
    and
    \begin{align}
        & |I_1| \lesssim \delta e^{(d+1)\Lambda}\left\|\rho\right\|_{L_t^\infty L_x^\infty} \lesssim \kappa, \label{eq:I1bound} \\
        & |I_2| \lesssim \bigg(1 + \frac{e^\Lambda}{\delta} + \frac{1}{\epsilon}\bigg)\left\|\rho_0\right\|_{W^{-1,1}} \lesssim \kappa, \label{eq:I2bound} \\
        & |I_3'| \lesssim T^\frac{d-1}{d} \frac{e^{d\Lambda}}{\epsilon} \left\|\rho\right\|_{L_t^\infty L_x^\infty} \left\|\nabla u\right\|_{L^1}\left(\left\|\phi*\varphi_{\delta'}\right\|_{L_t^2L_x^2}^2 - \left\|\phi*\varphi_{\delta'}*\varphi_{2\delta e^\Lambda}\right\|_{L_t^2L_x^2}^2 \right)^\frac{1}{d} \lesssim \kappa, \label{eq:I3'bound} \\
        & |I_4'| \lesssim Te^{2\Lambda}\left\|\rho\right\|_{L_t^\infty L_x^\infty} \left\|\nabla u\right\|_{L^1}^\frac{1}{2} \left(\left\|\nabla u*\varphi_{\epsilon'}\right\|_{L^1} - \left\|\nabla u * \varphi_{\epsilon'} * \varphi_{4\epsilon}\right\|_{L^1} \right)^\frac{1}{2} \lesssim \kappa, \label{eq:I4'bound} \\
        & |I_5'| \lesssim \frac{T}{\Lambda}\left\|\rho\right\|_{L_t^\infty L_x^\infty}\left\|\nabla u\right\|_{L^1} \lesssim \kappa, \label{eq:I5'bound} \\
        & |R_1| \lesssim T\frac{\delta'e^{(d+1)\Lambda}}{\epsilon\delta} \left\|\rho\right\|_{L_t^\infty L_x^\infty} \left\|\nabla u\right\|_{L^1} \lesssim \kappa, \label{eq:R1bound} \\
        & |R_2| \lesssim T \frac{\epsilon' e^{(d+1)\Lambda}}{\delta} \left\|\rho\right\|_{L_t^\infty L_x^\infty} \left\|\nabla u\right\|_{L^1} \lesssim \kappa, \label{eq:R2bound}
    \end{align}
    for some $\Lambda, \epsilon, \epsilon', \delta, \delta'>0$.

    We begin with the bound on $I_5'$ \eqref{eq:I5'bound}. This bound is satisfied if
    \begin{equation}\label{eq:Lambda}
        \Lambda=\frac{1}{\kappa}T\left\|\rho\right\|_{L_t^\infty L_x^\infty}\left\|\nabla u\right\|_{L^1},
    \end{equation}
    and will be the main dimensionless parameter. Indeed, we will see that it is convenient to rephrase most of the conditions on $\kappa$ in terms of $\Lambda$ instead.

    We now fix $\epsilon>0$ to be chosen later. Turning to the bound on $I_3'$ \eqref{eq:I3'bound}, the goal is to choose $\delta>0$ and $\delta'>0$ so that the following difference of norms is small:
    \begin{equation*}
        \left\|\phi*\varphi_{\delta'}\right\|_{L_t^2L_x^2}^2 - \left\|\phi*\varphi_{\delta'}*\varphi_{2\delta e^\Lambda}\right\|_{L_t^2L_x^2}^2,
    \end{equation*}
    while taking into account the bound on $R_1$ \eqref{eq:R1bound}, which requires the ratio $\frac{\delta'}{\delta}$ to be suitably small.
    
    To this end we fix some $\alpha>0$ to be chosen later, and define our \textit{pigeonholes} for $n \in \mathbb{N}$ as
    \begin{equation}\label{eq:pigeonhole1}
        \left\|\phi*\varphi_{\alpha^{n+1}}\right\|_{L_t^2L_x^2}^2 - \left\|\phi*\varphi_{\alpha^{n+1}}*\varphi_{\alpha^{n}}\right\|_{L_t^2L_x^2}^2 \ge 0,
    \end{equation}
    which are positive since $\varphi(h)\ge0$ is a positive mollifier \eqref{eq:mollifierassumptions}. The goal is now to control their sum. To this end, we wish to approximate them by the more summable differences
    \begin{equation*}
        \left\|\phi*\varphi_{\alpha^{n+1}}\right\|_{L_t^2L_x^2}^2 - \left\|\phi*\varphi_{\alpha^{n}}\right\|_{L_t^2L_x^2}^2.
    \end{equation*}
    The cross-term gives the error:
    \begin{align*}
        \left\|\phi*\varphi_{\alpha^{n}}\right\|_{L_t^2L_x^2}^2 - \left\|\phi*\varphi_{\alpha^{n+1}}*\varphi_{\alpha^{n}}\right\|_{L_t^2L_x^2}^2 & \lesssim \left\|\phi\right\|_{L_t^2L_x^2}\left(\left\|\phi*\varphi_{\alpha^{n}}-\phi*\varphi_{\alpha^{n+1}}*\varphi_{\alpha^{n}}\right\|_{L_t^2L_x^2}\right) \\
        & \lesssim \left\|\phi\right\|_{L_t^2L_x^2} \alpha^{n+1} \left\|\nabla(\phi*\varphi_{\alpha^{n}})\right\|_{L_t^2L_x^2} \\
        & \lesssim \left\|\phi\right\|_{L_t^2L_x^2}^2 \alpha^{n+1} \left\|\nabla\varphi_{\alpha^{n}}\right\|_{L^1} \\
        & \lesssim T \frac{\alpha^{n+1}}{\alpha^{n}} = T\alpha,
    \end{align*}
    where we have used the bound on the torus $\left\|\phi\right\|_{L_t^2L_x^2}^2\le T\left\|\phi\right\|_{L_t^\infty L_x^\infty}\le T$ \eqref{eq:phibound}, and the standard bound $\left\|\nabla\varphi_{\alpha^{n}}\right\|_{L^1}\lesssim\frac{1}{\alpha^{n}}$.

    We now bound the sum of the pigeonholes \eqref{eq:pigeonhole1}. For some $N_1\in\mathbb{N}$ to be chosen later:
    \begin{align*}
        \sum_{n=1}^{N_1} \left(\left\|\phi*\varphi_{\alpha^{n+1}}\right\|_{L_t^2L_x^2}^2 - \left\|\phi*\varphi_{\alpha^{n+1}}*\varphi_{\alpha^{n}}\right\|_{L_t^2L_x^2}^2\right) & \lesssim \sum_{n=1}^{N_1} \left(\left\|\phi*\varphi_{\alpha^{n+1}}\right\|_{L_t^2L_x^2}^2 - \left\|\phi*\varphi_{\alpha^{n}}\right\|_{L_t^2L_x^2}^2 + T\alpha\right) \\
        & = \left\|\phi*\varphi_{\alpha^{N_1+1}}\right\|_{L_t^2L_x^2}^2 - \left\|\phi*\varphi_{\alpha}\right\|_{L_t^2L_x^2}^2 + TN_1\alpha \\
        & \lesssim T(1 + N_1\alpha),
    \end{align*}
    where in the last line, we have again used the bound $\left\|\phi\right\|_{L_t^2L_x^2}^2\le T\left\|\phi\right\|_{L_t^\infty L_x^\infty}\le T$ \eqref{eq:phibound}.

    Since each pigeonhole \eqref{eq:pigeonhole1} is positive, there exists some
    \begin{equation}\label{eq:n1}
        n_1\in\{1,\dots,N_1\},
    \end{equation}
    such that
    \begin{equation*}
        \left\|\phi*\varphi_{\alpha^{n_1+1}}\right\|_{L_t^2L_x^2}^2 - \left\|\phi*\varphi_{\alpha^{n_1+1}}*\varphi_{\alpha^{n_1}}\right\|_{L_t^2L_x^2}^2 \lesssim T\left(\frac{1}{N_1}+\alpha\right).
    \end{equation*}

    Taking now
    \begin{align}
        & \delta = \frac{1}{2}e^{-\Lambda}\alpha^{n_1}, \label{eq:delta} \\
        & \delta' = \alpha^{n_1+1}, \label{eq:delta'}
    \end{align}
    then the bound on $I_3'$ \eqref{eq:I3'bound} is satisfied if
    \[
    T^\frac{d-1}{d} \frac{e^{d\Lambda}}{\epsilon} \left\|\rho\right\|_{L_t^\infty L_x^\infty} \left\|\nabla u\right\|_{L^1}  T^{\frac{1}{d}}\left(\frac{1}{N_1}+\alpha\right)^{\frac{1}{d}}\lesssim \kappa,
    \]
    that is 
    \[
    \frac{\Lambda e^{d\Lambda }}{\epsilon} \left(\frac{1}{N_1}+\alpha\right)^{\frac{1}{d}}\lesssim 1 \Leftrightarrow \frac{1}{N_1}+\alpha\lesssim 
    \epsilon^d \Lambda^{-d} e^{-d^2 \Lambda}
    \]
    for example if
    \begin{align}
        & \alpha = \epsilon^d e^{-d(d+1)\Lambda}, \label{eq:alpha} \\
        & N_1 = \left\lceil\epsilon^{-d} e^{d(d+1)\Lambda}\right\rceil, \label{eq:N1}
    \end{align}
    where we have used that $e^{-d\Lambda}\lesssim\Lambda^{-d}$.

    To summarise, given fixed $\epsilon>0$, then we have chosen $\alpha>0$ \eqref{eq:alpha}, $N_1\in\mathbb{N}$ \eqref{eq:N1}, then $n_1\in\{1,\dots,N_1\}$ \eqref{eq:n1}, and then $\delta>0$ \eqref{eq:delta}, $\delta'>0$ \eqref{eq:delta'} so that the bound on $I_3'$ \eqref{eq:I3'bound} is satisfied.

   % We now impose the condition \begin{align}\label{eq:epsiloncondition}\epsilon \le 1, T\left\|\nabla u\right\|_{L^1}.\end{align}

    We now show that $\delta \le \epsilon e^{-\Lambda}$ \eqref{eq:deltabound} using the expression for $\delta>0$ \eqref{eq:delta}, this is satisfied if
    \begin{align*}
        & \frac{1}{2}\alpha^{n_1}\le \epsilon \Longleftrightarrow \frac{1}{2} \epsilon^{n_1 d} e^{- d(d+1)n_1\Lambda}\le \epsilon,
    \end{align*}
    and this is satisfied if $\epsilon\le1$.
Next we prove that the bounds on $I_1$ \eqref{eq:I1bound}, and on $R_1$ \eqref{eq:R1bound} hold true if
    \begin{equation}\label{eq:epsiloncondition}
        \epsilon \le \min\left(1,\frac{\kappa}{\left\|\rho\right\|_{L_t^\infty L_x^\infty}}\right).
    \end{equation}
    First for $I_1$ \eqref{eq:I1bound}:
    \[
    \delta e^{(d+1)\Lambda}\left\|\rho\right\|_{L_t^\infty L_x^\infty} \lesssim \kappa,
    \]
    using the expression for $\delta>0$ \eqref{eq:delta}, this is equivalent to
    \[
     \frac{1}{2}e^{-\Lambda}\alpha^{n_1} e^{(d+1)\Lambda} \lesssim \frac{\kappa}{\left\|\rho\right\|_{L_t^\infty L_x^\infty}},
  \]
   hence from the expression for $\alpha>0$ \eqref{eq:alpha} we get 
   \[
   \epsilon^{n_1 d} e^{-n_1 d (d+1) \Lambda}  e^{d\Lambda} \lesssim \frac{\kappa}{\left\|\rho\right\|_{L_t^\infty L_x^\infty}},
   \]
   and this is satisfied if $\epsilon\le \min\Big(1,\frac{\kappa}{\left\|\rho\right\|_{L_t^\infty L_x^\infty}}\Big)$.
    
    We turn now to the bound on $R_1$ \eqref{eq:R1bound}. Using the expression for $\Lambda>0$ \eqref{eq:Lambda} we get 
    \begin{equation*}
        T\frac{\delta'e^{(d+1)\Lambda}}{\epsilon\delta} \left\|\rho\right\|_{L_t^\infty L_x^\infty} \left\|\nabla u\right\|_{L^1} \lesssim \kappa \Longleftrightarrow
        \delta'\lesssim \epsilon\delta \Lambda^{-1}e^{-(d+1)\Lambda} \Longleftrightarrow \alpha^{n_1+1} \lesssim \epsilon \alpha^{n_1} e^{-\Lambda}\Lambda^{-1}e^{-(d+1)\Lambda},
    \end{equation*}
   where we used the expression for $\delta>0$ \eqref{eq:delta}, $\delta'>0$ \eqref{eq:delta'}, hence using the expression for $\alpha>0$ \eqref{eq:alpha}
    \begin{equation*}
        \epsilon^d e^{-d(d+1)\Lambda} \lesssim \epsilon \Lambda^{-1} e^{-(d+2)\Lambda}.
    \end{equation*}
    and this is satisfied if $\epsilon \le 1$.

    % To summarise so far, given fixed $\epsilon>0$ satisfying
    % \begin{equation}\label{eq:epsiloncondition2}
    %     \epsilon\le \min(1,T\left\|\nabla u\right\|_{L^1}),
    % \end{equation}
    % we have chosen $\delta>0$ \eqref{eq:delta}, $\delta'>0$ \eqref{eq:delta'} so that the condition $\delta\le\epsilon e^{-\Lambda}$ \eqref{eq:deltabound}, and the bounds on $I_1$ \eqref{eq:I1bound}, on $I_3'$ \eqref{eq:I3'bound}, and on $R_1$ \eqref{eq:R1bound} are satisfied.

    Before finally turning to the bounds on $I_4'$ \eqref{eq:I4'bound} and $R_2$ \eqref{eq:R2bound}, we first bound the size of $\delta>0$ \eqref{eq:delta} in terms of the remaining parameter $\epsilon>0$. Using the expression for $\alpha>0$ \eqref{eq:alpha} and the requirement $\epsilon \le 1$:
    \begin{align*}
        \delta  \le \alpha  \le \epsilon.
    \end{align*}
    Next, we give a lower bound in terms of $\epsilon>0$. To this end, it is useful to define
    \begin{equation*}
        x = \frac{1}{\epsilon} e^{(d+1)\Lambda},
    \end{equation*}
    and then by the expression for $\delta>0$ \eqref{eq:delta}, $\alpha>0$ \eqref{eq:alpha}, and $N_1\in\mathbb{N}$ \eqref{eq:N1}:
    \begin{align*}
        \frac{1}{\delta} & \lesssim e^\Lambda \left(\frac{1}{\alpha}\right)^{N_1}= e^\Lambda \left(\frac{1}{\epsilon^de^{-d(d+1)\Lambda}}\right)^{\lceil\epsilon^{-d}e^{d(d+1)\Lambda}\rceil} \lesssim e^\Lambda (x^d)^{x^d+1} \lesssim e^\Lambda \exp(\exp(x)),
    \end{align*}
    % where we have used that $(x^d)^{x^d+1}\lesssim \exp(\exp(x))$ for $x\ge0$, 
    hence
    \begin{equation}\label{eq:deltasize}
        \frac{1}{\epsilon} \le \frac{1}{\delta} \lesssim e^{\Lambda} \exp\left(\exp\left(\frac{1}{\epsilon} e^{(d+1)\Lambda}\right)\right).
    \end{equation}

    Turning now to the bound on $I_4'$ \eqref{eq:I4'bound}, the goal is to choose $\epsilon>0$ and $\epsilon'>0$ so that the following difference of norms is small:
    \begin{equation}\label{eq:pigeonhole2bound}
        \left\|\nabla u * \varphi_{\epsilon'} \right\|_{L^1} - \left\|\nabla u*\varphi_{\epsilon'}*\varphi_{4\epsilon}\right\|_{L^1} \lesssim \Lambda^{-2}e^{-4\Lambda} \left\|\nabla u\right\|_{L^1},
    \end{equation}
    with the ratio $\frac{\epsilon'}{\delta}$ small to satisfy the bound on $R_2$ \eqref{eq:R2bound}, for which we will use the (double) exponential bound on $\frac{1}{\delta}$ \eqref{eq:deltasize}. Inspired by this, we fix some $x_0>0, M\in\mathbb{N}$ to be chosen later, and define our \textit{pigeonholes} for $n\in\mathbb{N}$ as
    \begin{equation}\label{eq:pigeonhole2}
        \left\|\nabla u * \varphi_{\left(\exp^{M(n+1)}(x_0)\right)^{-1}} \right\|_{L^1} - \left\|\nabla u*\varphi_{\left(\exp^{M(n+1)}(x_0)\right)^{-1}}*\varphi_{\left(\exp^{Mn}(x_0)\right)^{-1}}\right\|_{L^1} \ge 0,
    \end{equation}
    where $\exp^n(x)$ refers to repeated exponentiation $\exp(\exp(\dots(\exp(x))))$ (tetration). Each difference \eqref{eq:pigeonhole2} is positive since $\varphi(h)\ge0$ is a positive mollifier \eqref{eq:mollifierassumptions}. The goal is now to control their sum.

    To this end, we wish to approximate them by the more summable differences
    \begin{equation*}
        \left\|\nabla u * \varphi_{\left(\exp^{M(n+1)}(x_0)\right)^{-1}} \right\|_{L^1} - \left\|\nabla u * \varphi_{\left(\exp^{Mn}(x_0)\right)^{-1}}\right\|_{L^1}.
    \end{equation*}
    
    The cross-term gives the error:
    \begin{align*}
        & \left\|\nabla u*\varphi_{\left(\exp^{Mn}(x_0)\right)^{-1}}\right\|_{L^1} - \left\|\nabla u*\varphi_{\left(\exp^{M(n+1)}(x_0)\right)^{-1}}*\varphi_{\left(\exp^{Mn}(x_0)\right)^{-1}}\right\|_{L^1} \\
        & \qquad \lesssim \frac{1}{\exp^{M(n+1)}(x_0)}\left\|\nabla u*\nabla\varphi_{\left(\exp^{Mn}(x_0)\right)^{-1}}\right\|_{L^1} \\
        & \qquad \lesssim \frac{1}{\exp^M\left(\exp^{Mn}(x_0)\right)}\left\|\nabla u\right\|_{L^1} \left\|\nabla\varphi_{\left(\exp^{Mn}(x_0)\right)^{-1}}\right\|_{L^1} \\
        & \qquad \lesssim \frac{\exp^{Mn}(x_0)}{\exp^M\left(\exp^{Mn}(x_0)\right)}\left\|\nabla u\right\|_{L^1} \\
        & \qquad \lesssim \frac{1}{M}\left\|\nabla u\right\|_{L^1},
    \end{align*}
    where we have used the crude bound $Nx\le \exp^N(x)$ for any $N\in\mathbb{N}$. We will see that this does not lose us much.

    We now bound the sum of the pigeonholes \eqref{eq:pigeonhole2}. For some $N_2\in\mathbb{N}$ to be chosen later:
    \begin{align*}
        & \sum_{n=1}^{N_2} \left(\left\|\nabla u * \varphi_{\left(\exp^{M(n+1)}(x_0)\right)^{-1}} \right\|_{L^1} - \left\|\nabla u*\varphi_{\left(\exp^{M(n+1)}(x_0)\right)^{-1}}*\varphi_{\left(\exp^{Mn}(x_0)\right)^{-1}}\right\|_{L^1}\right) \\
        & \qquad \lesssim \sum_{n=1}^{N_2} \left(\left\|\nabla u * \varphi_{\left(\exp^{M(n+1)}(x_0)\right)^{-1}} \right\|_{L^1} - \left\|\nabla u * \varphi_{\left(\exp^{Mn}(x_0)\right)^{-1}}\right\|_{L^1} + \frac{1}{M}\left\|\nabla u\right\|_{L^1}\right) \\
        & \qquad = \left\|\nabla u * \varphi_{\left(\exp^{M(N_2+1)}(x_0)\right)^{-1}} \right\|_{L^1} - \left\|\nabla u * \varphi_{\left(\exp^{M}(x_0)\right)^{-1}}\right\|_{L^1} + \frac{N_2}{M}\left\|\nabla u\right\|_{L^1} \\
        & \qquad \lesssim \left(1 + \frac{N_2}{M}\right)\left\|\nabla u\right\|_{L^1}.
    \end{align*}

    Since each pigeonhole \eqref{eq:pigeonhole2} is positive, there exists some
    \begin{equation}\label{eq:n2}
        n_2\in\{1,\dots,N_2\},
    \end{equation}
    such that
    \begin{equation*}
        \left\|\nabla u * \varphi_{\left(\exp^{M(n_2+1)}(x_0)\right)^{-1}} \right\|_{L^1} - \left\|\nabla u*\varphi_{\left(\exp^{M(n_2+1)}(x_0)\right)^{-1}}*\varphi_{\left(\exp^{Mn_2}(x_0)\right)^{-1}}\right\|_{L^1} \lesssim \left(\frac{1}{N_2} + \frac{1}{M}\right)\left\|\nabla u\right\|_{L^1}.
    \end{equation*}

    Taking now
    \begin{align}
        & \epsilon = \frac{1}{4}\left(\exp^{Mn_2}(x_0)\right)^{-1}, \label{eq:epsilon} \\
        & \epsilon' = \left(\exp^{M(n_2+1)}(x_0)\right)^{-1}, \label{eq:epsilon'}
    \end{align}
    then the bound on $I_4'$ \eqref{eq:I4'bound} is satisfied if
    \[
    Te^{2\Lambda}\left\|\rho\right\|_{L_t^\infty L_x^\infty}\left\|\nabla u\right\|_{L^1}\left(\frac{1}{N_2} + \frac{1}{M}\right)^\frac{1}{2} \lesssim \kappa.
    \]
    Equivalently, using the expression for $\Lambda>0$ \eqref{eq:Lambda} this becomes
    \[
    \Lambda e^{2\Lambda} \left(\frac{1}{N_2} + \frac{1}{M}\right)^\frac{1}{2} \lesssim 1 \Leftrightarrow \frac{1}{N_2} + \frac{1}{M} \lesssim 
    \Lambda^{-2} e^{-4 \Lambda}.
    \]
    This bound is satisfied if
    \begin{align}
        & M,N_2 \ge e^{5\Lambda}, \label{eq:M,N2}
    \end{align}
    where we have used that $\Lambda^2 \lesssim e^{\Lambda}$.

    It remains to ensure the bound on $R_2$ \eqref{eq:R2bound} is satisfied, namely:
    \begin{equation}\label{eq:deltatetration}
        T \frac{\epsilon' e^{(d+1)\Lambda}}{\delta} \left\|\rho\right\|_{L_t^\infty L_x^\infty} \left\|\nabla u\right\|_{L^1} \lesssim \kappa \iff \frac{1}{\delta} \lesssim \frac{1}{\epsilon'}\Lambda^{-1} e^{(d+1)\Lambda} \iff \frac{1}{\delta} \lesssim \exp^{M(n_2+1)}(x_0) \Lambda^{-1} e^{(d+1)\Lambda},
    \end{equation}
    where we used the expression for $\epsilon'>0$ \eqref{eq:epsilon'}. Hence, using the bound for $\delta>0$ \eqref{eq:deltasize} it is sufficient for
    \begin{equation*}
        e^\Lambda \exp\left(\exp\left(\frac{1}{\epsilon}e^{(d+1)\Lambda}\right)\right) \lesssim \exp^{M(n_2+1)}(x_0) \Lambda^{-1} e^{(d+1)\Lambda}.
    \end{equation*}
    Using now the expression for $\epsilon>0$ \eqref{eq:epsilon}, and again the crude bound $Nx \le \exp^N(x)$ for any $N\in\mathbb{N}$, this becomes
    \begin{align*}
        & e^\Lambda \exp\left(\exp\left(e^{(d+1)\Lambda}\exp^{Mn_2}(x_0)\right)\right) \lesssim \exp^{M(n_2+1)}(x_0) \Lambda^{-1} e^{(d+1)\Lambda} \\
        & \impliedby \exp^{ \left\lceil e^{(d+1)\Lambda}\right\rceil + Mn_2}(x_0) \lesssim \exp^{M(n_2+1)}(x_0) \\
        & \impliedby  \left\lceil e^{(d+1)\Lambda}\right\rceil + Mn_2 \le M(n_2+1) ,
        %\\ & \impliedby  e^{(d+1)\Lambda} \le M,
    \end{align*}
    % for example if
    % \begin{equation*}
    %     M \ge 4e^{(d+1)\Lambda}.
    % \end{equation*}
    hence in order to satisfy also \eqref{eq:M,N2}, we therefore take
    \begin{align}
        & N_2 = \left\lceil e^{5\Lambda} \right\rceil, \label{eq:N2} \\
        & M = \left\lceil  e^{(d+5)\Lambda} \right\rceil. \label{eq:M}
    \end{align}

    It remains to pick $x_0>0$ to satisfy the smallness condition on $\epsilon>0$ \eqref{eq:epsiloncondition}:
    \begin{equation*}
        \epsilon \le \min\left(1,\frac{\kappa}{\left\|\rho\right\|_{L_t^\infty L_x^\infty}}\right).
    \end{equation*}
    Using the expression for $\epsilon>0$ \eqref{eq:epsilon} this is satisfied if
    \begin{equation}\label{eq:x0}
        x_0 = \frac{1}{\kappa}\left\|\rho\right\|_{L_t^\infty L_x^\infty}.
    \end{equation}

    Therefore, we are left only to satisfy the bound on $I_2$ \eqref{eq:I2bound}, which we now rephrase as a smallness condition on the weak norm $\left\|\rho_0\right\|_{W^{-1,1}}$, namely:
    \begin{equation*}
        \left\|\rho_0\right\|_{W^{-1,1}} \lesssim \kappa \left(1+\frac{e^\Lambda}{\delta}+\frac{1}{\epsilon}\right)^{-1}.
    \end{equation*}
    Using the bound on the size of $\delta\le\epsilon\le1$ \eqref{eq:deltasize}, \eqref{eq:epsiloncondition}, it is sufficient for
    \begin{equation*}
        \left\|\rho_0\right\|_{W^{-1,1}} \lesssim \kappa \delta e^{-\Lambda}.
    \end{equation*}

    Using that we have satisfied the bound \eqref{eq:deltatetration} on the size of $\frac{1}{\delta}$, it is even sufficient for
    \begin{equation*}
         \left\|\rho_0\right\|_{W^{-1,1}} \lesssim \kappa \left(\exp^{M(n_2+1)}(x_0)\right)^{-1}\Lambda e^{-(d+2)\Lambda},
    \end{equation*}
    and using $n_2\le N_2$, and the expressions for $M\in\mathbb{N}$ \eqref{eq:M}, $N_2\in\mathbb{N}$ \eqref{eq:N2}, $x_0>0$ \eqref{eq:x0}:
    \begin{align*}
        & \left\|\rho_0\right\|_{W^{-1,1}} \lesssim \kappa \left(\exp^{M(n_2+1)}(x_0)\right)^{-1}\Lambda e^{-(d+2)\Lambda} \\
        & \impliedby \left\|\rho_0\right\|_{W^{-1,1}} \lesssim \kappa \Lambda e^{-(d+2)\Lambda} \left(\exp^{M(N_2+1)}(x_0)\right)^{-1} \\
        & \iff \left\|\rho_0\right\|_{W^{-1,1}} \lesssim \kappa \Lambda e^{-(d+2)\Lambda} \left(\exp^{\left\lceil e^{(d+5)\Lambda}\right\rceil\left(\left\lceil e^{5\Lambda}\right\rceil+1\right)}\left(\frac{1}{\kappa}\left\|\rho\right\|_{L_t^\infty L_x^\infty}\right)\right)^{-1}.
    \end{align*}
    
    Finally, we denote by $C>0$ the constant (depending only on the dimension $d\ge2$ and the previously chosen mollifier \eqref{eq:mollifierassumptions}) such that under this assumption then $\left\|\rho_T\right\|_{W^{-1,1}}\le C\kappa$ \eqref{eq:kappabound}. Then by the crude bound $Nx\le\exp^N(x)$ for $N\in\mathbb{N}$, and the expression for $\Lambda>0$ \eqref{eq:Lambda}, it is sufficient for
    \begin{equation*}
        \left\|\rho_0\right\|_{W^{-1,1}} \le C\kappa \left(\exp^{\left\lceil A\exp\left(\frac{B}{\kappa}\left\|\rho\right\|_{L_t^\infty L_x^\infty}T\left\|\nabla u\right\|_{L^1}\right)\right\rceil}\left(\frac{1}{\kappa}\left\|\rho\right\|_{L_t^\infty L_x^\infty}\right)\right)^{-1}
    \end{equation*}
    for some constants $A,B>0$ depending only on the dimension $d\ge2$ (and the previously chosen of mollifier $\varphi(h)\in C_c^\infty(\mathbb{R}^d;\mathbb{R})$ satisfying \eqref{eq:mollifierassumptions}). Redefining now $\kappa \mapsto \frac{\kappa}{C}$ and $B \mapsto CB$ then gives the result.

    The proof of Corollary \ref{thm:mixingBV}, on the other hand, follows by the contrapositive of Theorem \ref{thm:quantitativeBV}, and reversing time so that $\rho_0(x)$ becomes $\rho_T(x)$ and vice versa. Namely, if
    \begin{equation*}
        \left\|\rho_0\right\|_{W^{-1,1}} > \kappa,
    \end{equation*}
    then
    \begin{equation*}
        \left\|\rho_T\right\|_{W^{-1,1}} > \kappa \left(\exp^{\left\lceil A\exp\left(\frac{B}{\kappa}\left\|\rho_0\right\|_{L^\infty}T\left\|\nabla u\right\|_{L^1}\right)\right\rceil}\left(\frac{1}{\kappa}\left\|\rho_0\right\|_{ L^\infty}\right)\right)^{-1},
    \end{equation*}
    Taking now $\kappa = \left\|\rho_0\right\|_{W^{-1,1}}-\alpha$ and letting $\alpha\to0$ then gives the result (after redefining $A\mapsto A+1$ since the ceiling function $\lceil\cdot\rceil$ is not upper continuous).
\end{proof}

\section{Appendix}

\subsection{Proof of Proposition \ref{prop:nonquantitativeBV}}

\begin{proof}

    For the sake of contradiction suppose that there exists some $M \ge 0$, $\epsilon>0$, and a sequence
    \begin{align*}
        & \left\{\rho_0^{(n)}\right\}_{n\in\mathbb{N}}\subset L^\infty(\mathbb{T}^d;\mathbb{R}), \\
        & \left\{u^{(n)}\right\}_{n\in\mathbb{N}}\subset BV(\mathbb{T}^d;\mathbb{R}^d),
    \end{align*}
    such that, denoting by $\rho^{(n)}(x,t) \in L_t^\infty L_x^\infty$ the \textit{unique} bounded solution (see \cite{ambrosio2004transport}) to the transport equation along $u^{(n)}(x)$ with initial datum $\rho_0^{(n)}(x)$, we have
    \begin{align}
        & \left\|\rho_0^{(n)}\right\|_{L^\infty} \le M \text{ for all } n \in\mathbb{N}, \nonumber \\
        & \left\|\nabla u^{(n)}\right\|_{L_x^1} \le M \text{ for all } n \in\mathbb{N}, \nonumber \\
        & \left\|\rho_0^{(n)}\right\|_{W^{-1,1}}\xrightarrow{n\to\infty}0, \label{eqmixhyp} \\
        & \left\|\rho^{(n)}_T\right\|_{W^{-1,1}} > \epsilon. \label{eqmixhyp2}
    \end{align}

    By Sobolev embedding, since $u^{(n)}(x) \in BV$ then also $u^{(n)}(x)\in L_x^\frac{d}{d-1}$ are uniformly bounded, and $\rho^{(n)}(x,t) \in L_t^\infty L_x^\infty$ are uniformly bounded, so also $\frac{\partial \rho^{(n)}}{\partial t}(x,t) \in L_t^\infty W_x^{-1,\frac{d}{d-1}}$ are uniformly bounded by the transport equation; $\frac{\partial \rho^{(n)}}{\partial t}(x,t)=-\nabla\cdot(u^{(n)}(x,t)\rho^{(n)}(x,t))$. The goal is then to apply the Aubin-Lions compactness lemma \cite{simon1986compact}, to show that the sequence $\rho^{(n)}(x,t)$ is compact in $L_t^\infty C_x^{-1}$ where $C^{-1}$ is the dual of $BV$. Working on the compact domain $\mathbb{T}^d$, we have the continuous embeddings
    \begin{equation*}
        L^\infty \subset C^{-1} \subset W^{-1,\frac{d}{d-1}},
    \end{equation*}
    where $C^{-1}$ is the Banach space dual of $BV$, and in particular the embedding $L^\infty \subset C^{-1}$ is compact. Then by the Aubin-Lions compactness lemma, the sequence $\rho^{(n)}(x,t)$ is compactly embedded in $C_t^0 C_x^{-1}$.
    
    Therefore, we see that $\rho^{(n)}(x,t)$ has a subsequence converging strongly in $C_t^0 C_x^{-1}$, with the limit also in $L_t^\infty L_x^\infty$. Without loss of generality, let this subsequence be the original sequence $\rho^{(n)}(x,t)$. Denote the limit by $\bar{\rho}(x,t) \in C_t^0 C_x^{-1} \cap L_t^\infty L_x^\infty$. Likewise, taking a further subsequence if necessary, we see that $ u^{(n)}(x)\in BV$ converges weakly-$*$ in $BV$ to some limit $\bar{u}(x)\in BV$.

    Since $C^{-1}$ is the dual of $BV$, it then follows that also $u^{(n)}(x) \rho^{(n)}(x,t) \xrightharpoonup{n\to\infty} \bar{u}(x)\bar{\rho}(x,t)$ converges in the sense of distribution (i.e. against smooth test functions) to $\bar{u}(x)\bar{\rho}(x,t) \in L_t^\infty L_x^1$. Note also that $\rho_0^{(n)}(x)\xrightharpoonup{n\to\infty}0$ converges to zero, strongly in $W^{-1,1}$, by \eqref{eqmixhyp}.
    
    Then, for any test function $\phi(x,t) \in C_c^\infty(\mathbb{T}^d\times[0,T);\mathbb{R})$ we have
    \begin{align*}
        & \int_{\mathbb{T}^d\times[0,T)} \bar{\rho}(x,t) \left(\frac{\partial \phi}{\partial t}(x,t) + u(x,t)\cdot\nabla \phi(x,t)\right) \; dxdt \\
        & \qquad = \lim_{n\to\infty} \int_{\mathbb{T}^d\times[0,T)} \rho^{(n)}(x,t) \frac{\partial \phi}{\partial t}(x,t) + \rho^{(n)}(x,t)u^{(n)}(x,t)\cdot\nabla \phi(x,t) \; dxdt \\
        & \qquad = - \lim_{n\to\infty} \int_{\mathbb{T}^d} \rho_0^{(n)}(x) \phi(x,0) \;dx, \\
        & \qquad = 0,
    \end{align*}
    and so $\bar{\rho}(x,t)\in L_t^\infty L_x^\infty$ is the \textit{unique} bounded solution \cite{ambrosio2004transport} to the transport equation along $\bar{u}(x) \in BV$ with zero initial datum, that is $\bar{\rho}(x,t)=0$ and in particular $\bar{\rho}(x,T)=0$.
    
    However, since $\rho^{(n)}_T(x)$ converges to $\bar{\rho}_T(x)$ strongly in $W^{-1,1}$ (since $\rho^{(n)}(x,t)$ converges strongly in $C_t^0 C_x^{-1}$), and $\left\|\rho^{(n)}_T\right\|_{W^{-1,1}} > \epsilon$ uniformly by \eqref{eqmixhyp2}, then also
    \begin{equation*}
        \left\|\bar{\rho}_T\right\|_{W^{-1,1}} \ge \epsilon>0,
    \end{equation*}
    a contradiction as required.
\end{proof}

\subsection{Proof of Corollary \ref{cor: proof of Ambrosio}}\label{proof: proof of Ambrosio}

\begin{proof}
    We use Proposition \ref{prop:fiveterms}, that is we show for any fixed test function $\phi_T(x) \in W^{1,\infty}$ that
    \begin{equation*}
        \int_{\mathbb{T}^d} \rho_T(x)\phi_T(x)\;dx=0,
    \end{equation*}
    and note that we may effectively take $T>0$ arbitrarily.

    To this end, we choose parameters as in Proposition \ref{prop:fiveterms} such that each of the bounds is arbitrarily small. Fix some arbitrarily small parameter $\kappa>0$.

    We first choose $\Lambda>0$ large enough that the bound on $I_5$ is less than or equal to $\kappa$, for instance
    \begin{equation*}
        \Lambda = \kappa T \left\|\rho\right\|_{L_t^\infty L_x^\infty} \left\|\nabla u\right\|_{L^1}.
    \end{equation*}
    Next, we pick $\epsilon >0$ small enough so that $I_4$ is also bounded by $\kappa$. By Proposition \ref{prop:fiveterms}, it is sufficient to satisfy
    \begin{equation*}
        Te^{2\Lambda}\left\|\rho\right\|_{L_t^\infty L_x^\infty} \left\|\nabla u\right\|_{L^1}^\frac{1}{2} \left(\left\|\nabla u\right\|_{L^1} - \left\|\nabla u * \varphi_{4\epsilon}\right\|_{L^1} \right)^\frac{1}{2} \le \kappa,
    \end{equation*}
    where now $\Lambda>0$ is fixed. Thus, we need only show that as $\epsilon\to 0$ that
    \begin{equation*}
        \left\|\nabla u * \varphi_{4\epsilon}\right\|_{L^1} \to \left\|\nabla u\right\|_{L^1}.
    \end{equation*}
    First, since $\varphi(x)\ge0$ is a positive mollifier, we have the bound
    \begin{equation*}
        \left\|\nabla u * \varphi_{4\epsilon}\right\|_{L^1} \le \left\|\nabla u\right\|_{L^1},
    \end{equation*}
    and so it is sufficient to prove
    \begin{equation*}
        \liminf \left\|\nabla u * \varphi_{4\epsilon}\right\|_{L^1} \ge \left\|\nabla u\right\|_{L^1},
    \end{equation*}
    which follows from convergence in the distributional sense of
    \begin{equation*}
        (\nabla u * \varphi_{4\epsilon})(x) \rightharpoonup \nabla u(x),
    \end{equation*}
    as $\epsilon\to0$.

    We have so far chosen $\Lambda>0$ large and $\epsilon>0$ small so that
    \begin{align*}
        & |I_5| \lesssim \kappa, \\
        & |I_4| \lesssim \kappa.
    \end{align*}

    By assumption $\rho_0(x)=0$ and so by Proposition \ref{prop:fiveterms}, $I_2=0$. We are left to pick $\delta>0$ small so that the assumption $\delta \le \epsilon e^{-\Lambda}$ of Proposition \ref{prop:fiveterms} is satisfied, and
    \begin{align*}
        & |I_1| \lesssim \kappa, \\
        & |I_3| \lesssim \kappa.
    \end{align*}
    Therefore, we require $\delta>0$ small enough that
    \begin{align*}
        & \delta \le \epsilon e^{-\Lambda}, \\
        & \delta e^{(d+1)\Lambda}\left\|\rho\right\|_{L_t^\infty L_x^\infty} \le \kappa, \\
        & T^\frac{d-1}{d} \frac{e^{d\Lambda}}{\epsilon} \left\|\rho\right\|_{L_t^\infty L_x^\infty} \left\|\nabla u\right\|_{L^1}\left(\left\|\phi\right\|_{L_t^2L_x^2}^2 - \left\|\phi*\varphi_{2\delta e^\Lambda}\right\|_{L_t^2L_x^2}^2 \right)^\frac{1}{d} \le \kappa,
    \end{align*}
    where only convergence of the last constraint requires further analysis. Since $\phi(x,t) \in L_t^\infty L_x^\infty \subset L_t^2 L_x^2$ (on the domain $\mathbb{T}^d\times[0,T]$), this follows by the standard convergence of
    \begin{equation*}
        (\phi*\varphi_{2\delta e^\Lambda})(x,t) \to \phi(x,t),
    \end{equation*}
    in $L_t^2L_x^2$ as $2\delta e^\Lambda \to 0$.

    Thus we may chose $\Lambda>0$ large, $\epsilon>0$ small, and $\delta>0$ small so that by Proposition \ref{prop:fiveterms}:
    \begin{align*}
        & |I_1| \lesssim \kappa, \\
        & |I_2| = 0, \\
        & |I_3| \lesssim \kappa, \\
        & |I_4| \lesssim \kappa, \\
        & |I_5| \lesssim \kappa.
    \end{align*}
    Since $\kappa>0$ is arbitrary, we deduce that
    \begin{align*}
        \int_{\mathbb{T}^d} \rho_T(x)\phi_T(x) \;dx & = I_1 + I_2 + I_3 + I_4 + I_5 \\
        & = 0.
    \end{align*}
    Since also $\phi_T(x) \in W^{1,\infty}$ is arbitrary, we deduce that
    \begin{equation*}
        \rho_T(x) = 0.
    \end{equation*}
    Since also $T>0$ is essentially arbitrary, we deduce that
    \begin{equation*}
        \rho(x,t)=0.
    \end{equation*}
\end{proof}

\addcontentsline{toc}{section}{References}
\bibliographystyle{apa}
\bibliography{references.bib}

% \bib, bibdiv, biblist are defined by the amsrefs package.
\begin{bibdiv}
\begin{biblist}

\bib{ambrosio2004transport}{article}{
      author={Ambrosio, Luigi},
       title={Transport equation and {C}auchy problem for {$BV$} vector fields},
        date={2004},
        ISSN={0020-9910,1432-1297},
     journal={Invent. Math.},
      volume={158},
      number={2},
       pages={227\ndash 260},
         url={https://doi.org/10.1007/s00222-004-0367-2},
      review={\MR{2096794}},
}

\bib{ambrosio2014continuity}{article}{
      author={Ambrosio, Luigi},
      author={Crippa, Gianluca},
       title={Continuity equations and {ODE} flows with non-smooth velocity},
        date={2014},
        ISSN={0308-2105},
     journal={Proc. Roy. Soc. Edinburgh Sect. A},
      volume={144},
      number={6},
       pages={1191\ndash 1244},
         url={https://doi.org/10.1017/S0308210513000085},
      review={\MR{3283066}},
}

\bib{ambrosio2005traces}{article}{
      author={Ambrosio, Luigi},
      author={Crippa, Gianluca},
      author={Maniglia, Stefania},
       title={Traces and fine properties of a {$BD$} class of vector fields and applications},
        date={2005},
        ISSN={0240-2963,2258-7519},
     journal={Ann. Fac. Sci. Toulouse Math. (6)},
      volume={14},
      number={4},
       pages={527\ndash 561},
         url={http://afst.cedram.org/item?id=AFST_2005_6_14_4_527_0},
      review={\MR{2188582}},
}

\bib{bianchini2020uniqueness}{article}{
      author={Bianchini, Stefano},
      author={Bonicatto, Paolo},
       title={A uniqueness result for the decomposition of vector fields in {$\mathbb{R}^d$}},
        date={2020},
        ISSN={0020-9910,1432-1297},
     journal={Invent. Math.},
      volume={220},
      number={1},
       pages={255\ndash 393},
         url={https://doi.org/10.1007/s00222-019-00928-8},
      review={\MR{4071413}},
}

\bib{bianchini2022differentiability}{article}{
      author={Bianchini, Stefano},
      author={De~Nitti, Nicola},
       title={Differentiability in measure of the flow associated with a nearly incompressible {BV} vector field},
        date={2022},
        ISSN={0003-9527,1432-0673},
     journal={Arch. Ration. Mech. Anal.},
      volume={246},
      number={2-3},
       pages={659\ndash 734},
         url={https://doi.org/10.1007/s00205-022-01820-1},
      review={\MR{4514062}},
}

\bib{bressan2003lemma}{article}{
      author={Bressan, Alberto},
       title={A lemma and a conjecture on the cost of rearrangements},
        date={2003},
        ISSN={0041-8994},
     journal={Rend. Sem. Mat. Univ. Padova},
      volume={110},
       pages={97\ndash 102},
      review={\MR{2033002}},
}

\bib{cooperman2023exponential}{article}{
      author={Cooperman, William},
       title={Exponential mixing by shear flows},
        date={2023},
        ISSN={0036-1410,1095-7154},
     journal={SIAM J. Math. Anal.},
      volume={55},
      number={6},
       pages={7513\ndash 7528},
         url={https://doi.org/10.1137/22M1513861},
      review={\MR{4665036}},
}

\bib{crippa2008estimates}{article}{
      author={Crippa, Gianluca},
      author={De~Lellis, Camillo},
       title={Estimates and regularity results for the {D}i{P}erna-{L}ions flow},
        date={2008},
        ISSN={0075-4102,1435-5345},
     journal={J. Reine Angew. Math.},
      volume={616},
       pages={15\ndash 46},
         url={https://doi.org/10.1515/CRELLE.2008.016},
      review={\MR{2369485}},
}

\bib{diperna1989ordinary}{article}{
      author={DiPerna, R.~J.},
      author={Lions, P.-L.},
       title={Ordinary differential equations, transport theory and {S}obolev spaces},
        date={1989},
        ISSN={0020-9910},
     journal={Invent. Math.},
      volume={98},
      number={3},
       pages={511\ndash 547},
         url={https://doi.org/10.1007/BF01393835},
      review={\MR{1022305}},
}

\bib{evans2018measure}{book}{
      author={Evans, Lawrence~C.},
      author={Gariepy, Ronald~F.},
       title={Measure theory and fine properties of functions},
     edition={Revised},
      series={Textbooks in Mathematics},
   publisher={CRC Press, Boca Raton, FL},
        date={2015},
        ISBN={978-1-4822-4238-6},
      review={\MR{3409135}},
}

\bib{hadvzic2018singular}{article}{
      author={Had\v{z}i\'{c}, Mahir},
      author={Seeger, Andreas},
      author={Smart, Charles~K.},
      author={Street, Brian},
       title={Singular integrals and a problem on mixing flows},
        date={2018},
        ISSN={0294-1449,1873-1430},
     journal={Ann. Inst. H. Poincar\'{e} C Anal. Non Lin\'{e}aire},
      volume={35},
      number={4},
       pages={921\ndash 943},
         url={https://doi.org/10.1016/j.anihpc.2017.09.001},
      review={\MR{3795021}},
}

\bib{huysmans_2024}{thesis}{
      author={Huysmans, Lucas},
       title={Passive scalar transport by non-smooth incompressible fluids: Mixing and vanishing viscosity},
        type={Ph.D. Thesis},
        date={2024},
         url={https://www.repository.cam.ac.uk/handle/1810/377378},
}

\bib{huysmans2024quantitative}{article}{
      author={Huysmans, Lucas},
      author={Said, Ayman~Rimah},
       title={Quantitative {E}stimates in {P}assive {S}calar {T}ransport: {A} {U}nified {A}pproach in {$W^{1,p}$} via {C}hrist-{J}ourné {C}ommutator {E}stimates},
        date={2025},
     journal={arXiv preprint arXiv:2402.11642},
}

\bib{iyer2014lower}{article}{
      author={Iyer, Gautam},
      author={Kiselev, Alexander},
      author={Xu, Xiaoqian},
       title={Lower bounds on the mix norm of passive scalars advected by incompressible enstrophy-constrained flows},
        date={2014},
        ISSN={0951-7715,1361-6544},
     journal={Nonlinearity},
      volume={27},
      number={5},
       pages={973\ndash 985},
         url={https://doi.org/10.1088/0951-7715/27/5/973},
      review={\MR{3207161}},
}

\bib{leger2018new}{article}{
      author={L\'{e}ger, Flavien},
       title={A new approach to bounds on mixing},
        date={2018},
        ISSN={0218-2025,1793-6314},
     journal={Math. Models Methods Appl. Sci.},
      volume={28},
      number={5},
       pages={829\ndash 849},
         url={https://doi.org/10.1142/S0218202518500215},
      review={\MR{3799259}},
}

\bib{seis2013maximal}{article}{
      author={Seis, Christian},
       title={Maximal mixing by incompressible fluid flows},
        date={2013},
        ISSN={0951-7715,1361-6544},
     journal={Nonlinearity},
      volume={26},
      number={12},
       pages={3279\ndash 3289},
         url={https://doi.org/10.1088/0951-7715/26/12/3279},
      review={\MR{3141856}},
}

\bib{simon1986compact}{article}{
      author={Simon, Jacques},
       title={Compact sets in the space {$L^p(0,T;B)$}},
        date={1987},
        ISSN={0003-4622},
     journal={Ann. Mat. Pura Appl. (4)},
      volume={146},
       pages={65\ndash 96},
         url={https://doi.org/10.1007/BF01762360},
      review={\MR{916688}},
}

\end{biblist}
\end{bibdiv}
\end{document}